\PassOptionsToPackage{sort&compress}{NJDnatbib}
\documentclass[ASNA,onecolumn]{USG}
\usepackage{mathrsfs,booktabs,enumitem,courier,needspace}
\usepackage[utf8]{inputenc}
\usepackage{numcompress}
\usepackage{adjustbox}
\usepackage{array}

\usepackage[section]{placeins}
\usepackage{etoolbox}
\AtBeginEnvironment{thebibliography}{\interlinepenalty=10000}
\usepackage[capitalize,nameinlink]{cleveref}
\usetikzlibrary{arrows.meta,positioning,fit,backgrounds,calc}
\hypersetup{colorlinks=true,citecolor=blue,linkcolor=blue,urlcolor=blue}
\floatstyle{ruled}
\restylefloat{algorithm}
\floatname{algorithm}{Algorithm}

\makeatletter
\let\c@theorem\relax

\makeatother
\theoremstyle{WBstyleone}
\newtheorem{theorem}{Theorem}[section]
\newtheorem{lemma}[theorem]{Lemma}
\newtheorem{proposition}[theorem]{Proposition}
\newtheorem{corollary}[theorem]{Corollary}
\theoremstyle{definition}
\newtheorem{definition}[theorem]{Definition}
\newtheorem{assumption}[theorem]{Assumption}
\newtheorem{remark}[theorem]{Remark}
\crefname{assumption}{Assumption}{Assumptions}

\numberwithin{equation}{section}

\newcommand{\Cc}{\mathcal{C}}

\newcommand{\Pp}{\mathscr{P}}
\newcommand{\dd}{\,\mathrm{d}}

\hypersetup{pdfauthor={Hayman Thabet},pdftitle={Conservative Implicit and Explicit Solvers for Reactive Conformable Transport with Multirate Memory}}
\RemoveFromHook{shipout/after}[USG]
\makeatletter
\renewcommand\appendix{\par
  \setcounter{section}{0}\setcounter{subsection}{0}\setcounter{equation}{0}%
  \gdef\thefigure{\@Alph\c@section.\arabic{figure}}%
  \gdef\thetable{\@Alph\c@section.\arabic{table}}%
  \gdef\thesection{\appendixname~\@Alph\c@section}%
  \gdef\theequation{\@Alph\c@section.\arabic{equation}}%
  \gdef\theHsection{appendix.\Alph{section}}}
\historydates{}
\renewcommand{\@@address}[2][]{%
  \g@addto@macro\jmkaddress{\hsize=\textwidth
    \stepcounter{affn}\xdef\@currentlabel{\theaffn}\jmkLabel{#1}%
    \addressfont #2\par}}
\renewcommand{\@maketitle}{%
  {\hsize=\textwidth\parindent=0pt
   \thispagestyle{titlepage}%
   \let\footnote\thanks
   \ifx\@title\empty\else{\@title\par}\fi
   \removelastskip\vskip8pt
   {\artauthors\par}%
   \removelastskip\vskip27.25pt
   \printabstract
   \removelastskip\vskip\adjtitleskip}}
\def\oddfoot@titlepage@info{\head@headings@def{\hfill\oddpage@headings@def}}
\def\evenfoot@titlepage@info{\oddfoot@titlepage@info}
\def\evenhead@headings@info{\head@headings@def{\evenpage@headings@def\hfill}}
\def\oddhead@titlepage@info{}
\def\evenhead@titlepage@info{}
\def\@@MSC[#1]#2{\def\@msccode{{\FIfont\textbf{#1 MSC: }#2\par}}}
\apptocmd{\printabstractpart}{\ifx\@msccode\@empty\else\@msccode\fi}{}{\PackageError{Manuscript}{MSC layout patch failed}{}}
\makeatother
\date{}
\usepackage{float}
\begin{document}

\title{Conservative Implicit and Explicit Solvers for Reactive
	Conformable Transport with Multirate Memory}

\author[1]{Zahraa Abdallah}
\author[2]{Hayman Thabet}

\address[1]{Independent Researcher, Beirut, Lebanon}
\address[2]{School of Mathematics and Data Science,
	University of New England, Biddeford, ME 04005, USA}

\corres{\email{haymanthabet@gmail.com}}

\keywords{Causal memory  |  residual stability  |  mobile--immobile transport  |  finite volumes  |  compartment approximation  |  mass conservation}
\MSC[2020]{65M12, 65M08, 35R09}

\abstract[ABSTRACT]{
In this article, we analyze conservative implicit and explicit schemes for reactive
transport with a conformable time scale and multirate exchange.
The conformable derivative defines a local change of time coordinate,
while exchange with immobile regions supplies the memory.
Positive Gauss--Jacobi quadrature approximates the prescribed exchange
spectrum by finite compartments. Both schemes use conservative
finite-volume fluxes. Forward Euler gives the explicit update;
backward Euler initializes the implicit second-order backward
differentiation formula (BDF2).
Eliminating the compartment states reduces each implicit step to a
tridiagonal system. Under stated mesh and regularity conditions, we
prove implicit solvability, conditional explicit positivity and
stability, and discrete mass balance.
A cell-average error estimate separates spatial, compartment and
temporal contributions. Its proof accounts for the boundary fluxes
and the BDF2 startup error. Source-specific bounds make the temporal
and memory prefactors uniform under mesh and compartment refinement.
For a broader positive-history model, continuous residual estimates
have coefficients independent of memory strength and kernel mass.
Local coercivity gives a sharper coefficient.
Manufactured tests recover second-order spatial convergence,
first-order Euler convergence and second-order BDF2 convergence.
Separate experiments examine the memory and combined estimates.
In a prescribed porous-medium application, the finest numerical
states differ by $3.57\times10^{-5}$ in relative space--time $L^2$.
Independent resolution checks support the computed pulse attenuation,
immobile storage and release after injection ends.
}

\maketitle

\section{Introduction}\label{sec:intro}
Transport through heterogeneous porous media often involves retention
and release over several time scales. Mobile solute exchanges mass
with immobile regions, whose distribution of exchange rates produces
a hereditary response. Multirate mass-transfer models describe this
mechanism \cite{HaggertyGorelick1995}.
Eliminated history formulations and finite-volume implementations
are established \cite{WangZhengGorelick2005,MunicchiEtAl2021}.
Silva et al.\ \cite{SilvaEtAl2009} derive analytical updates for
immobile concentrations in a real-time formulation.
Wang et al.\ \cite{WangCarreraEtAl2022} use block elimination in
nonlinear reactive transport.
These approaches provide the modeling and computational background
for the present study.

Our focus is the error analysis of two conservative schemes for a
specified exchange spectrum. Spatial resolution, compartment
quadrature and time stepping introduce different errors.
Their effects must be separated to interpret refinement results.
The boundary fluxes also require attention: local boundary defects
can persist while cell-average state errors converge at second order.
Mass conservation supplies another check, but does not measure
state accuracy.

The analysis begins with a broader model that combines a local
conformable time scale with a positive history operator.
For differentiable states, the conformable derivative is local
\cite{khalil2014conformable,abdeljawad2015conformable}.
The coordinate $\vartheta=H(t)$ in Eq.~\eqref{eq:conformable}
converts it to an ordinary derivative.
The separate integral operator supplies memory
\cite{tarasov2018nononlocality,ortigueira2015whatis}.
Maximum principles, energy estimates and representations for related
nonlocal evolution equations appear in
\cite{luchko2009maximum,vergara2015decay,du2017nonlocal,du2020stochastic}.
Here we use the sign of the history term at a running maximum of
the error norm. This yields residual estimates without integrating
the full kernel mass at the initial endpoint.

For the transport application, a positive capacity density determines
the exchange law. Gauss--Jacobi quadrature replaces this density by
finite compartments. Gaussian discretizations of diffusive memory
\cite{MonteghettiMatignonPiot2020} and exponential-kernel compression
\cite{Baffet2019} provide related approximation tools.
We combine the compartment rule with centered finite-volume fluxes
and either forward Euler or backward-Euler-started BDF2.
The spatial and multistep ingredients are classical
\cite{EymardGallouetHerbin2000,Dahlquist1978}.
The results below concern their stability and error propagation for
the specified memory law and boundary treatment.

\Needspace{5\baselineskip}
The contributions are as follows.
\begin{enumerate}[label=(\arabic*),leftmargin=*,itemsep=3pt]
\item \textbf{Continuous residual estimates.}
We bound state error by a weighted continuous residual.
The coefficients are independent of memory strength and kernel mass.
Local coercivity gives a sharper coefficient. For zero initial error,
it is optimal over the stated operator class
(Theorems~\ref{thm:stability} and~\ref{thm:coercive}).

\item \textbf{Discrete solver properties.}
We prove solvability of the eliminated implicit system.
Under the stated mesh condition, a computable step restriction
ensures explicit positivity and weighted stability. Both schemes satisfy a discrete mass identity
when cumulative quantities use the state update's time rule
(Propositions~\ref{prop:p3pair-implicit-solvability}--\ref{prop:p3pair-balance}).

\item \textbf{A combined cell-average error estimate.}
The estimate separates spatial, compartment and temporal errors.
An elliptic projection treats the actual boundary fluxes.
An integrated kernel defect controls compartment error.
The temporal estimate retains the backward-Euler startup term
for BDF2 (Theorem~\ref{thm:p3ce-combined}).

\item \textbf{Uniform source bounds and separate verification.}
For the prescribed manufactured and pulse sources, we bound the
temporal and memory prefactors independently of mesh and compartment
resolution. The explicit admissible step still depends on stiffness.
Independent implementation checks and separate refinement studies
test the predicted behavior. Application outputs receive their own
resolution checks.
\end{enumerate}

The application uses prescribed nondimensional coefficients and
injection data. It examines concentration, detector response and
mass inventories without calibration to field measurements.
Section~\ref{sec:model} defines the model and assumptions.
Section~\ref{sec:method} presents the schemes, and
Section~\ref{sec:analysis} gives their analysis.
Section~\ref{sec:results} reports the numerical studies.
Section~\ref{sec:conclusion} discusses the findings and their limits.

The main notation is summarized in Table~\ref{tab:p3-principal-notation}.
\begin{table}[!htbp]
\centering
\small
\caption{Main notation. All model quantities are nondimensional.}
\label{tab:p3-principal-notation}
\setlength{\tabcolsep}{5pt}
\renewcommand{\arraystretch}{1.15}
\begin{tabular}{@{}>{\raggedright\arraybackslash}p{0.27\linewidth}>{\raggedright\arraybackslash}p{\dimexpr0.73\linewidth-10pt\relax}@{}}
\toprule
Symbol & Definition \\
\midrule
$t,\ \vartheta,\ H$ & Model time, physical time, and clock map: $\vartheta=H(t)$. \\
$\alpha(t),\ \psi(t)$ & Local temporal order and multiplier $\psi(t)=t^{1-\alpha(t)}$. \\
$u,\ z_\lambda$ & Mobile concentration and immobile concentration at exchange rate $\lambda$. \\
$v,\ D(x)$ & Transport velocity and dispersion coefficient. \\
$\kappa_d,\ \kappa_m$ & Mobile degradation rate and trapping amplitude. \\
$q(x,\vartheta)$ & Prescribed source of mobile solute. \\
$b(\lambda),\ b_j,\ \lambda_j$ & Continuous capacity density, compartment capacities, and compartment exchange rates. \\
$\ell,\ \ell_M$ & Continuous relaxation kernel and its $M$-compartment approximation. \\
$N_x,\ h$ & Number of spatial cells and cell width $h=1/N_x$. \\
$M=N_c$ & Number of immobile compartments. \\
$N_\vartheta,\ \tau,\ \Theta$ & Time-step count, physical step $\tau=\Theta/N_\vartheta$, and final physical time ($\Theta=\vartheta_f$ in the application). \\
$U_h,\ U_{h,M},\ U_{h,M}^n$ & Mobile states on the same mesh: continuous time and memory; continuous time with $M$ compartments; and fully discrete at $\vartheta_n=n\tau$, respectively. \\
$L_h,\ \Pi_h$ & Conservative finite-volume transport operator and exact cell-averaging operator. \\
$\|V\|_{2,h}$ & Discrete spatial norm, $(h\sum_i|V_i|^2)^{1/2}$. \\
$D_M,\ \varepsilon_M^{(1)}$ & Primitive kernel discrepancy $D_M(\vartheta)=\int_0^\vartheta(\ell-\ell_M)(s)\,ds$ and its $L^1(0,\Theta)$ norm. \\
$R^c[w]$ & Continuous equation residual of a trial state $w$. \\
$M_m,\ M_i$ & Total mobile and capacity-weighted immobile mass. \\
$M_{\rm cap},\ M_{\rm deg}$ & Cumulative boundary capture and degradation loss. \\
\bottomrule
\end{tabular}
\end{table}

\section{Model and assumptions}\label{sec:model}
Let $\Omega\subset\mathbb R^d$ be a bounded Lipschitz domain.
Set $Q_T=\Omega\times(0,T]$, with $T>0$.
Space and time are nondimensional.
Unless specified otherwise, inner products and norms are those of
$L^2(\Omega)$. For $\alpha\in C([0,T])$ satisfying
$0<\alpha_{\min}\le\alpha(t)\le1$, define
\begin{equation}\label{eq:conformable}
 \psi(t)=t^{1-\alpha(t)},\qquad T_t^{\alpha(t)}v=\psi(t)v_t,
 \qquad H(t)=\int_0^t\psi(r)^{-1}\dd r.
\end{equation}
For $0<t\le\min\{1,T\}$, $\psi(t)^{-1}=t^{\alpha(t)-1}\le t^{\alpha_{\min}-1}$. Since $\alpha_{\min}>0$, this bound is integrable at zero. Away from zero, $\psi^{-1}$ is continuous. Thus $H(T)<\infty$. Since $H'(t)=\psi(t)^{-1}>0$, the coordinate $\vartheta=H(t)$ transforms $\psi v_t$ into $\mathrm dv/\mathrm d\vartheta$. This identity holds for variable and constant local order.

Let $\gamma$ be piecewise continuous, right-continuous at zero, and satisfy $0\le\gamma(t)\le\gamma_{\max}<1$. The prescribed modulation is chosen from
\begin{equation}\label{eq:admissible}
 \Pp=\left\{\rho\in C([0,T]\times[0,1]):\rho\ge0,\quad
 \int_0^1\rho(t,s)\dd s=1\right\}.
\end{equation}
Write $\rho_{\max}=\max_{[0,T]\times[0,1]}\rho$, which is finite by continuity. Its associated kernel and causal-difference operator are
\begin{align}
 \omega_\rho(t,\tau)&=\frac{(t-\tau)^{-\gamma(t)}}{t\,\Gamma(1-\gamma(t))}
 \rho(t,\tau/t), \quad 0\le\tau<t,\label{eq:kernel}\\
 \Cc_\rho[v](x,t)&=\int_0^t\omega_\rho(t,\tau)
[v(x,t)-v(x,\tau)]\dd\tau.\label{eq:correction}
\end{align}
The mass of the full kernel obeys
\begin{equation}\label{eq:kernel_mass}
 0\le k_\rho(t):=\int_0^t\omega_\rho(t,\tau)\dd\tau
 \le \frac{\rho_{\max}t^{-\gamma(t)}}{\Gamma(2-\gamma(t))}.
\end{equation}
The normalization in \eqref{eq:admissible} fixes the integral of the modulation $\rho$ over relative time. The full mass $k_\rho(t)$ also depends on the endpoint factor and can grow as $t\downarrow0$. It is finite at every positive time. The stability argument uses pointwise finiteness and nonnegativity. It does not require integrability of $\psi^{-1}k_\rho$.

We consider the initial--boundary value problem
\begin{equation}\label{eq:model}
 \mathcal L_\rho[u]:=\psi u_t+c\cdot\nabla u-\nu\Delta u
             +\eta\Cc_\rho[u]+\mu\phi(u)=s,
 \qquad u(\cdot,0)=u_0,\quad u|_{\partial\Omega}=g.
\end{equation}
Here $c$ is the advection field and $\nu\ge\nu_0>0$ is the diffusion coefficient. The coefficient $\eta\ge0$ sets the memory strength. The reaction is $\mu\phi$. The source is denoted by $s(x,t)$. In history integrals, the dummy variable $s\in[0,1]$ denotes relative time. Examples of normalized modulations are $\rho\equiv1$, $\rho(t,s)=2(1-s)$ and, for $\lambda>0$,
\[
 \rho_\lambda(t,s)=\frac{\lambda t e^{-\lambda t(1-s)}}{1-e^{-\lambda t}},
 \qquad \rho_\lambda(0,s)=1.
\]
The last choice introduces exponential tempering of the relative-time profile. Its use here specifies the modulation in \eqref{eq:kernel}, rather than a standard tempered fractional derivative \cite{sabzikar2015tempered}.

\subsection{Meaning of the memory term}
Normalization and positivity imply $k_\rho(t)>0$ for $t>0$. Introduce the weighted past state
\begin{equation}\label{eq:past_average}
 \overline u_\rho(t)=\int_0^t p_t(\tau)u(\tau)\dd\tau,
 \qquad p_t(\tau)=\omega_\rho(t,\tau)/k_\rho(t),\qquad
 \int_0^t p_t(\tau)\dd\tau=1.
\end{equation}
Then $\eta\Cc_\rho[u]=\eta k_\rho(t)(u-\overline u_\rho)$. In the evolution equation, this term drives the present state toward its weighted past. The modulation $\rho$ distributes weight over past times. The exponent $\gamma$ sets the endpoint singularity. The coefficient $\eta$ controls the memory strength. The factor $1/t$ expresses the modulation in the coordinate $\tau/t$. The full kernel mass still depends on time.

Under the local time change $t=H^{-1}(\vartheta)$, the memory kernel becomes
\[
 \widetilde\omega(\vartheta,\sigma)
 =\omega_\rho(t(\vartheta),t(\sigma))\psi(t(\sigma)),
 \qquad 0<\sigma<\vartheta.
\]
It remains causal and nonnegative, and generally depends on both transformed times. Thus changing the local clock does not eliminate the hereditary term. Setting $\eta=0$ removes memory, whereas $\alpha\equiv1$ restores the ordinary time derivative. These identities describe the model algebraically. Physical validation would require an independently justified constitutive law and measured observations.

\subsection{Regularity and reaction bounds}
The analysis compares an exact state with a trial state under the following sufficient assumptions.
\begin{assumption}\label{ass:standing}
The coefficients $c,\nu$ belong to $W^{1,\infty}(Q_T)$, $\nu\ge\nu_0>0$, and $\eta\ge0$. Let $K\subset\mathbb R$ contain the essential ranges of both states.
On $K$, the reaction satisfies
\begin{equation}\label{eq:onesided}
 \mu[\phi(a)-\phi(b)](a-b)\ge-\kappa_\phi(a-b)^2,
       \qquad a,b\in K,\quad \kappa_\phi\ge0.
\end{equation}
Assume that an exact solution exists in the following class.
Both compared states belong to
$C([0,T];L^2(\Omega))$. On every $[\varepsilon,T]$, with
$0<\varepsilon<T$, both are absolutely continuous in $L^2(\Omega)$.
They also belong to $L^2(\varepsilon,T;H^2(\Omega))$.
Their boundary traces agree. Their difference therefore belongs to
$H^1_0(\Omega)$ almost everywhere.
The reaction difference is locally integrable with values in
$L^2(\Omega)$. The error equation holds almost everywhere in that space.
Its forcing satisfies the weighted integrability required by each theorem.
These conditions justify the norm derivative and spatial integrations.
Exact initial data are imposed where stated.
\end{assumption}
In particular, existence is assumed. The conclusions below concern uniqueness and error control within this regularity class. For $\mu\ge0$, cubic reaction satisfies \eqref{eq:onesided} globally with $\kappa_\phi=0$. Quadratic reaction has the same constant on $[0,\infty)$; on $[m,\infty)$ with $m<0$, one can take $\kappa_\phi=2\mu|m|$. A range established only at sampled points does not supply the global set $K$ required by the estimate.

The next bound relates temporal increments to the size of the memory operator.
\begin{lemma}[History increment bound]\label{lem:props}
If $\|v(t)-v(\tau)\|\le L|t-\tau|^\sigma$, $0<\sigma\le1$, then
\begin{equation}\label{eq:history_holder}
 \|\Cc_\rho[v](t)\|\le
 \frac{\rho_{\max}L}{(1+\sigma-\gamma(t))\Gamma(1-\gamma(t))}
 t^{\sigma-\gamma(t)}.
\end{equation}
The operator is linear, causal, and annihilates time-independent states.
\end{lemma}
\begin{proof}
Fix $t>0$. The assumed increment bound makes the integrand Bochner integrable in $L^2(\Omega)$, because $\sigma-\gamma(t)>-1$. The triangle inequality for that integral gives
\begin{align*}
 \|\Cc_\rho[v](t)\|
 &\le \frac{1}{t\Gamma(1-\gamma(t))}
   \int_0^t (t-\tau)^{-\gamma(t)}\rho(t,\tau/t)
    \|v(t)-v(\tau)\|\dd\tau\\
 &\le \frac{\rho_{\max}L}{t\Gamma(1-\gamma(t))}
   \int_0^t(t-\tau)^{\sigma-\gamma(t)}\dd\tau\\
 &=\frac{\rho_{\max}L}{(1+\sigma-\gamma(t))\Gamma(1-\gamma(t))}
  t^{\sigma-\gamma(t)}.
\end{align*}
This proves \eqref{eq:history_holder}. For fixed $\rho$, the integral is linear in $v$. It uses only values at times $\tau\le t$, so it is causal. A time-independent state has zero increment at every history time and therefore has zero memory response.
\end{proof}

\subsection{Mobile--immobile transport application}
\label{sec:porous-application}
Let $x\in(0,1)$ and let $\vartheta$ denote nondimensional physical time. The
mobile concentration is $u(x,\vartheta)$, while $z_\lambda(x,\vartheta)$ is
the concentration in immobile regions with exchange rate
$\lambda>\lambda_0>0$. The nonnegative density $b(\lambda)$ specifies their
capacity relative to the mobile phase. The transport equations are
\begin{align}
	u_\vartheta+\partial_x j+\kappa_d u
	+\int_{\lambda_0}^{\infty}\lambda b(\lambda)
	(u-z_\lambda)\,d\lambda &= q(x,\vartheta),
	&j&=v u-D(x)u_x, \label{eq:app-mobile}\\
	\partial_\vartheta z_\lambda&=\lambda(u-z_\lambda).
	\label{eq:app-immobile}
\end{align}
Here $v>0$ is the transport velocity, $D(x)>0$ is the dispersion coefficient,
$\kappa_d\geq0$ is the mobile-phase degradation rate, and $q$ is a prescribed
injection source. Both phases are initially solute-free. The conditions
$u(0,\vartheta)=u(1,\vartheta)=0$ represent perfectly absorbing reservoirs
at the two faces. Consequently, a concentration trace at an interior detector
is distinguished from the boundary capture flux.

Choose $0<\gamma<1$ and a trapping amplitude $\kappa_m>0$, and set
\begin{equation}
	b(\lambda)=\frac{\kappa_m(\lambda-\lambda_0)^{\gamma-1}}
	{\Gamma(\gamma)\Gamma(1-\gamma)\lambda^2},
	\qquad \lambda>\lambda_0.
	\label{eq:app-capacity}
\end{equation}
Its total capacity is finite:
$\int_{\lambda_0}^{\infty}b(\lambda)\,d\lambda
=(1-\gamma)\kappa_m\lambda_0^{\gamma-2}$.
Eliminating the initially empty immobile regions gives
\begin{equation}
	\int_{\lambda_0}^{\infty}\lambda b(\lambda)(u-z_\lambda)\,d\lambda
	=\ell(\vartheta)u(x,\vartheta)
	+\int_0^\vartheta K(\vartheta-\sigma)
	[u(x,\vartheta)-u(x,\sigma)]\,d\sigma,
	\label{eq:app-elimination}
\end{equation}
where
\begin{equation}
	K(r)=\frac{\kappa_m e^{-\lambda_0 r}r^{-\gamma}}
	{\Gamma(1-\gamma)},\qquad
	\ell(\vartheta)=\int_\vartheta^\infty K(r)\,dr
	=\frac{\kappa_m\lambda_0^{\gamma-1}
		\Gamma(1-\gamma,\lambda_0\vartheta)}{\Gamma(1-\gamma)}.
	\label{eq:app-kernel}
\end{equation}
The upper incomplete gamma function is
$\Gamma(a,z)=\int_z^\infty r^{a-1}e^{-r}\,dr$.
The local term $\ell(\vartheta)u$ is part of the exchange law and must be
retained. In particular, the initially empty immobile regions are not
equivalent to regions initially equilibrated with a nonzero mobile state.

\subsubsection{Connection to the conformable model}
Use the model coordinate $t\in[0,1]$ and the prescribed clock
$\vartheta=H(t)=\int_0^t\psi(r)^{-1}\,dr$,
with $\alpha(t)=1-0.2t$ and $\psi(t)=t^{0.2t}$ for $t>0$;
define $\psi(0)=1$. Thus $e^{-0.2/e}\leq\psi\leq1$. Writing
$\widetilde u(x,t)=u(x,H(t))$ and suppressing the tilde below gives
$\partial_\vartheta u=\psi(t)\partial_tu$.
The same clock is used for both phases. It is a coordinate representation,
not an additional retention mechanism. Comparisons between clocks must
therefore use the same physical time $\vartheta$ and transform the source
and memory consistently.

For $t>0$ and $0\leq s<1$, define
\begin{align}
	\Delta_H(t,s)&=H(t)-H(ts),\nonumber\\
	B(t,s)&=\left[\frac{t(1-s)}{\Delta_H(t,s)}\right]^\gamma
	\frac{e^{-\lambda_0\Delta_H(t,s)}}{\psi(ts)},
	&Z(t)&=\int_0^1 B(t,s)\,ds,\nonumber\\
	\rho(t,s)&=B(t,s)/Z(t),
	&\eta(t)&=\kappa_m t Z(t).
	\label{eq:app-modulation}
\end{align}
The continuous endpoint extensions are $B(t,1)=\psi(t)^{\gamma-1}$
and $B(0,s)=1$. Hence $\rho\geq0$, $\int_0^1\rho\,ds=1$, and
$\eta\geq0$, with $\eta(0)=0$.
Writing $\mathcal C_\rho$ for the paper's causal-difference operator,
Eq.~\eqref{eq:app-mobile} becomes
\begin{equation}
	\psi u_t+(v-D')u_x-Du_{xx}
	+[\kappa_d+\ell(H(t))]u+\eta(t)\mathcal C_\rho[u]
	=q(x,H(t)).
	\label{eq:app-conformable}
\end{equation}
For the ordinary clock $H(t)=t$, Eq.~\eqref{eq:app-modulation} reduces to
\begin{equation}
	\rho(t,s)=\frac{\lambda_0t e^{-\lambda_0t(1-s)}}
	{1-e^{-\lambda_0t}},\qquad
	\eta(t)=\frac{\kappa_m}{\lambda_0}(1-e^{-\lambda_0t}).
	\label{eq:app-ordinary-clock}
\end{equation}
Thus the ordinary-clock specialization has a normalized exponential
modulation, with the coupling and local exchange term determined by
the capacity spectrum.

The conservative flux is essential when $D$ varies. To retain an admissible
spatial conformable representation, prescribe
\begin{equation}
	p(x)=\tfrac34[4x(1-x)]^3,\qquad
	\alpha_s(x)=2-p(x),\qquad
	D(x)=D_{\rm ref}x^{p(x)},\qquad D(0)=D_{\rm ref}.
	\label{eq:app-dispersion}
\end{equation}
Then $1.25\leq\alpha_s\leq2$, $D\in C^2([0,1])$, and
$D(x)\geq D_{\rm ref}e^{-256/243}>0$.
Consequently the coefficients $\nu=D$ and $c=v-D'$ satisfy the stated
spatial regularity requirements. The spatial order specifies a local
coefficient profile; it does not represent spatially nonlocal transport.

The continuous residual estimate for this conservative formulation is proved in Proposition~\ref{prop:app-residual-control}.

\subsubsection{Prescribed pulse and physical diagnostics}
For the baseline application take
\[
(v,D_{\rm ref},\kappa_d,\kappa_m,\lambda_0,\gamma)
=(2,0.02,0.1,0.8,0.5,0.5).
\]
Let $g(x)$ be $\exp[-((x-0.15)/0.04)^2]$ divided by its integral
over $(0,1)$. Prescribe $q(x,\vartheta)=g(x)P(\vartheta)$, where
\begin{equation}
	P(\vartheta)=
	\begin{cases}
		\displaystyle\frac{2}{0.20}
		\sin^2\!\left(\frac{\pi(\vartheta-0.05)}{0.20}\right),
		&0.05\leq\vartheta\leq0.25,\\
		0,&\text{otherwise}.
	\end{cases}
	\label{eq:app-pulse}
\end{equation}
The injected mass is one. The final physical time is $H(1)$.
The no-trapping case sets $\kappa_m=0$; in that case both the hereditary
term and $\ell$ vanish. It is a physical control, distinct from a comparison of numerical schemes. All other transport and injection parameters are
identical in the two physical cases.

The main observables are the interior trace $u(0.75,\vartheta)$, its
peak and peak time, the mobile mass $M_m=\int_0^1u\,dx$, and the
immobile mass $M_i=\int_0^1\int_{\lambda_0}^{\infty}b(\lambda)z_\lambda
\,d\lambda\,dx$. Define the cumulative captured mass and degradation loss by
\begin{equation}
	M_{\rm cap}(\vartheta)=\int_0^\vartheta[j(1,\sigma)-j(0,\sigma)]\,d\sigma,
	\qquad
	M_{\rm deg}(\vartheta)=\kappa_d\int_0^\vartheta M_m(\sigma)\,d\sigma.
	\label{eq:app-capture}
\end{equation}
The governing balance is
\begin{equation}
	M_m+M_i+M_{\rm cap}+M_{\rm deg}
	=\int_0^\vartheta P(\sigma)\,d\sigma.
	\label{eq:app-balance}
\end{equation}
Its discrepancy is a separate diagnostic. A small $L^2$ state error alone
does not bound the boundary flux, which involves a spatial derivative.
\section{Conservative implicit and explicit compartment schemes}
\label{sec:method}
The schemes in this section approximate the mobile--immobile system
\eqref{eq:app-mobile}--\eqref{eq:app-immobile}, with the capacity law
\eqref{eq:app-capacity}. They are not discretizations of every admissible
kernel in \eqref{eq:admissible}. Both use the same compartment quadrature,
spatial fluxes and physical time $\vartheta=H(t)$.

\subsection{Positive compartment quadrature}
\label{sec:p3pair-quadrature}
Under $r=\lambda_0/\lambda$, the capacity measure becomes
\[
 b(\lambda)\,d\lambda
 =C_b\frac{r^{1-\gamma}(1-r)^{\gamma-1}}
             {\mathrm B(2-\gamma,\gamma)}\,dr,\qquad
 C_b=(1-\gamma)\kappa_m\lambda_0^{\gamma-2},
\]
where the reversal of the integration limits is understood.
Let $(y_k,w_k)_{k=1}^{N_c}$ be Gauss--Jacobi nodes and weights for
$(1-y)^{\gamma-1}(1+y)^{1-\gamma}$ on $(-1,1)$. Define
\begin{equation}
 r_k=\frac{1+y_k}{2},\qquad
 \lambda_k=\frac{\lambda_0}{r_k},\qquad
 b_k=C_b\frac{w_k}{\sum_{j=1}^{N_c}w_j}.
 \label{eq:p3pair-rates}
\end{equation}
Then $\lambda_k>0$, $b_k>0$ for $\kappa_m>0$, and
$\sum_k b_k=C_b$. Positive Gaussian quadrature provides a finite
compartment approximation of the continuous rate distribution
\cite{gautschi2004orthogonal,MonteghettiMatignonPiot2020}.
The no-trapping control sets $\kappa_m=0$ and hence $b_k=0$.
Exact total capacity does not imply exact exchange dynamics: the factors
$\lambda_k$ and $e^{-\lambda_k\vartheta}$ must also be resolved.

\subsection{Conservative spatial discretization}
\label{sec:p3pair-space}
For $N_x\geq2$, let $h=1/N_x$, $x_i=(i-\tfrac12)h$ and
$x_{i+1/2}=ih$. The unknowns $U_i$ and $Z_{k,i}$ approximate mobile and
immobile cell averages. At interior faces set
\begin{equation}
 F_{i+1/2}(U)=\frac{v}{2}(U_i+U_{i+1})
       -\frac{D(ih)}{h}(U_{i+1}-U_i),\qquad 1\leq i<N_x.
 \label{eq:p3pair-flux}
\end{equation}
The absorbing boundary values give
\begin{equation}
 F_{1/2}(U)=-\frac{2D(0)}hU_1,\qquad
 F_{N_x+1/2}(U)=\frac{2D(1)}hU_{N_x},\qquad
 (L_hU)_i=-\frac{F_{i+1/2}(U)-F_{i-1/2}(U)}h.
 \label{eq:p3pair-boundary-flux}
\end{equation}
Let $Q_i(\vartheta)$ denote the cell average of $q(\cdot,\vartheta)$.
The semidiscrete equations are
\begin{align}
 \dot U&=L_hU-\kappa_dU-\sum_{k=1}^{N_c}b_k\lambda_k(U-Z_k)+Q,
 \label{eq:p3pair-semimobile}\\
 \dot Z_k&=\lambda_k(U-Z_k),\qquad k=1,\ldots,N_c.
 \label{eq:p3pair-semiimmobile}
\end{align}
The dot denotes a physical-time derivative. The application starts from
$U^0=Z_k^0=0$; the schemes also allow prescribed nonzero initial vectors.
For the positivity arguments below, impose
\begin{equation}
 |v|h\leq 2D(ih),\qquad 1\leq i<N_x.
 \label{eq:p3pair-mesh}
\end{equation}
This condition makes the off-diagonal entries of $L_h$ nonnegative.
Both methods retain this centered flux.

\subsection{Implicit scheme: backward-Euler startup and BDF2}
\label{sec:p3pair-implicit}
Use $\vartheta_n=n\tau$, $\tau=\vartheta_f/N_\vartheta$, and define
\begin{equation}
 \delta_nV^n=a_0V^n+a_1V^{n-1}+a_2V^{n-2},\qquad
 (a_0,a_1,a_2)=
 \begin{cases}
 (1,-1,0)/\tau,&n=1,\\
 (3/2,-2,1/2)/\tau,&n\geq2.
 \end{cases}
 \label{eq:p3pair-bdf}
\end{equation}
Apply this difference to both equations
\eqref{eq:p3pair-semimobile}--\eqref{eq:p3pair-semiimmobile}, evaluating
their right-hand sides at level $n$. Put
\begin{equation}
 B^n=-a_1U^{n-1}-a_2U^{n-2},\quad
 P_k^n=-a_1Z_k^{n-1}-a_2Z_k^{n-2},\quad
 \chi_n=\sum_k\frac{b_k\lambda_k a_0}{a_0+\lambda_k}.
 \label{eq:p3pair-known-history}
\end{equation}
Algebraic elimination of the compartment states gives
\begin{align}
 [(a_0+\kappa_d+\chi_n)I-L_h]U^n
 &=B^n+Q^n+\sum_k\frac{b_k\lambda_k}{a_0+\lambda_k}P_k^n,
 \label{eq:p3pair-implicit}\\
 Z_k^n&=\frac{\lambda_kU^n+P_k^n}{a_0+\lambda_k}.
 \label{eq:p3pair-recovery}
\end{align}
Terms at level $n-2$ are omitted at startup. The mobile solve is
tridiagonal, and recovery of all compartment states is algebraic.
This elimination is exact for the fully discrete compartment system.

\subsection{Explicit scheme: forward Euler}
\label{sec:p3pair-explicit}
On the same spatial and compartment discretization, use
\begin{align}
 U^{n+1}&=U^n+\tau\left[L_hU^n-\kappa_dU^n
       -\sum_k b_k\lambda_k(U^n-Z_k^n)+Q^n\right],
       \label{eq:p3pair-explicit}\\
 Z_k^{n+1}&=(1-\tau\lambda_k)Z_k^n+\tau\lambda_kU^n.
       \label{eq:p3pair-explicit-z}
\end{align}
Both updates use the old state, and the source is evaluated at
$\vartheta_n$. Put $S_{N_c}=\sum_k b_k\lambda_k$ and define
\begin{equation}
 \tau_{\rm pos}=
 \min\left\{\frac1{\max_k\lambda_k},
       \frac1{\max_i[\kappa_d+S_{N_c}-(L_h)_{ii}]}\right\}.
 \label{eq:p3pair-positive-step}
\end{equation}
The explicit implementation requires $0<\tau\leq\tau_{\rm pos}$.
For a requested safety factor $0<\sigma<1$, it selects
$N_\vartheta=\lceil\vartheta_f/(\sigma\tau_{\rm pos})\rceil$.
Section~\ref{sec:p3pair-analysis} proves the resulting positivity and
stability statements.

\subsection{Mass accounting and implementation cost}
\label{sec:p3pair-accounting}
Define
\begin{equation}
 M_m^n=h\sum_iU_i^n,\qquad
 M_i^n=h\sum_{i,k}b_kZ_{k,i}^n,\qquad
 \Phi^n=\frac2h[D(0)U_1^n+D(1)U_{N_x}^n].
 \label{eq:p3pair-masses}
\end{equation}
Initialize cumulative injection $I_{\rm inj}$, capture $M_{\rm cap}$ and
degradation $M_{\rm deg}$ at zero. For the implicit scheme use
\[
 \delta_n I_{\rm inj}^n=h\sum_i Q_i^n,\qquad
 \delta_n M_{\rm cap}^n=\Phi^n,\qquad
 \delta_n M_{\rm deg}^n=\kappa_d M_m^n.
\]
For forward Euler advance the same quantities with the same old-time
rates as the states. This synchronized accounting is required for the
balance identity proved below.

Each step costs $O(N_xN_c)$ operations; the implicit method also solves
one $O(N_x)$ tridiagonal system. Both methods require $O(N_xN_c)$ working
storage when only the necessary state levels are retained. Saved output
arrays add storage according to the chosen output frequency. No complete
time history is required during stepping. 

\section{Residual stability, discrete properties and convergence}
\label{sec:analysis}
We first derive continuous residual estimates and specialize them to
the conservative application. We then analyze the discrete schemes
and combine their spatial, compartment and temporal errors.
Each result states its required regularity and mesh conditions.
\subsection{A bound independent of the kernel mass}\label{sec:an_stab}
Let $u$ and $v$ satisfy the same boundary data and write $e=u-v$. Their difference is governed by
\begin{equation}\label{eq:error_equation}
 \psi e_t+c\cdot\nabla e-\nu\Delta e+\eta\Cc_\rho[e]
  +\mu[\phi(u)-\phi(v)]=f,
 \qquad e|_{\partial\Omega}=0.
\end{equation}
Define
\begin{equation}\label{eq:b}
 b=\tfrac12\|\nabla\!\cdot c\|_\infty
  +\frac{\|\nabla\nu\|_\infty^2}{2\nu_0}+\kappa_\phi.
\end{equation}
The coefficient norms in \eqref{eq:b} are essential suprema over $Q_T$, with spatial derivatives. The zero boundary trace gives
\[
 \langle c\cdot\nabla e,e\rangle
   =-\frac12\int_\Omega(\nabla\!\cdot c)e^2\dd x
   \ge-\frac12\|\nabla\!\cdot c\|_\infty\|e\|^2.
\]
For the nondivergence-form diffusion term, integration by parts and Young's inequality yield
\begin{align*}
 -\langle\nu\Delta e,e\rangle
 &=\int_\Omega\nu|\nabla e|^2\dd x
   +\int_\Omega e\nabla\nu\cdot\nabla e\dd x\\
 &\ge\nu_0\|\nabla e\|^2
   -\|\nabla\nu\|_\infty\|e\|\|\nabla e\|\\
 &\ge\frac{\nu_0}{2}\|\nabla e\|^2
       -\frac{\|\nabla\nu\|_\infty^2}{2\nu_0}\|e\|^2.
\end{align*}
After adding the integrated one-sided reaction bound, the local operator is bounded below by $\nu_0\|\nabla e\|^2/2-b\|e\|^2$. At a running maximum of $\|e(t)\|$, the entire causal-difference term has a favorable sign. Retaining that sign yields the following estimate.

\begin{theorem}[Running-maximum stability]\label{thm:stability}
Let Assumption~\ref{ass:standing} hold and let
$\psi^{-1}\|f(\cdot,t)\|\in L^1(0,T)$. Then
\begin{equation}\label{eq:stability_l1}
 \max_{0\le t\le T}\|e(t)\|
 \le e^{bH(T)}\|e(0)\|
       +\int_0^T e^{b[H(T)-H(r)]}\psi(r)^{-1}\|f(r)\|\dd r.
\end{equation}
In particular, if $\psi^{-1/2}f\in L^2(Q_T)$,
\begin{equation}\label{eq:stability}
 \max_t\|e(t)\|\le e^{bH(T)}\|e(0)\|
+C_R\|\psi^{-1/2}f\|_{L^2(Q_T)},\qquad
 C_R=\begin{cases}
 \sqrt{(e^{2bH(T)}-1)/(2b)},&b>0,\\
 \sqrt{H(T)},&b=0.
 \end{cases}
\end{equation}
For fixed $b$ and clock $H$, the displayed coefficients are independent
of $\eta$ and $k_\rho$. No integrability of $\psi^{-1}k_\rho$ near zero
is required.
\end{theorem}
\begin{proof}
Set $y(t)=\|e(t)\|$ and $Y(t)=\max_{0\le r\le t}y(r)$.
Both are continuous on $[0,T]$, and $y$ is absolutely continuous on
$[\varepsilon,T]$ for every $0<\varepsilon<T$. For
$\varepsilon\le s<t\le T$,
\[
 0\le Y(t)-Y(s)\le\int_s^t|y'(r)|\dd r.
\]
Indeed, if $Y(t)>Y(s)$, choose $r_*\in(s,t]$ with $y(r_*)=Y(t)$
and use $Y(s)\ge y(s)$. Applying this bound to disjoint intervals
shows that $Y$ is also absolutely continuous on $[\varepsilon,T]$.

For each $t>0$, the history integral is Bochner integrable because
\[
 \int_0^t\omega_\rho(t,r)\|e(t)-e(r)\|\dd r
 \le 2k_\rho(t)Y(t)<\infty.
\]
At a contact time $y(t)=Y(t)$, positivity of the kernel gives
\[
 \langle\Cc_\rho[e](t),e(t)\rangle
 \ge y(t)\int_0^t\omega_\rho(t,r)[y(t)-y(r)]\dd r\ge0.
\]
At almost every such time with $y(t)>0$,
$\langle e_t,e\rangle=yy'$. Testing
\eqref{eq:error_equation} and using the local lower bound yields
$\psi yy'\le by^2+\|f\|y$.

On $\{y<Y\}$, $Y$ is locally constant. On $\{y=Y\}$, the
nonnegative absolutely continuous function $Y-y$ has derivative zero
almost everywhere, so $Y'=y'$ there. Also, $Y'=0$ almost everywhere
on $\{Y=0\}$. Division by $y$ only where $y>0$, together with $b\ge0$,
therefore gives
\[
 Y'\le\psi^{-1}(bY+\|f\|)
 \quad\text{almost everywhere on }[\varepsilon,T].
\]
Multiplication by $e^{-bH}$ and integration imply
\[
 Y(T)\le e^{b[H(T)-H(\varepsilon)]}Y(\varepsilon)
 +\int_\varepsilon^T
 e^{b[H(T)-H(r)]}\psi(r)^{-1}\|f(r)\|\dd r.
\]
As $\varepsilon\downarrow0$, continuity gives
$Y(\varepsilon)\to\|e(0)\|$ and $H(\varepsilon)\to0$.
The integrand is bounded by
$e^{bH(T)}\psi^{-1}\|f\|\in L^1(0,T)$, so dominated convergence
proves \eqref{eq:stability_l1}.

If $\psi^{-1/2}f\in L^2(Q_T)$, Cauchy--Schwarz and $H(T)<\infty$
give the required weighted $L^1$ integrability. Applying
Cauchy--Schwarz to the forcing term in \eqref{eq:stability_l1}
gives \eqref{eq:stability}, since
\[
 \int_0^T e^{2b[H(T)-H(r)]}\psi(r)^{-1}\dd r
 =\int_0^{H(T)}e^{2b[H(T)-h]}\dd h=C_R^2.
\]
Only the sign and fixed-time finiteness of the history kernel were used.
\end{proof}

\begin{corollary}[Uniqueness and residual control]\label{cor:error}
For fixed coefficients, data and kernel, there is at most one solution
in the class and common range of Assumption~\ref{ass:standing}.
Let $u$ be an exact solution and let $w$ be a trial state such that
$(u,w)$ satisfies the regularity and range conditions of that assumption.
Suppose that their initial and boundary data agree. Define
$R^c[w]=\mathcal L_\rho[w]-s$. If $\psi^{-1/2}R^c[w]\in L^2(Q_T)$, then
\begin{equation}
 \max_{0\leq t\leq T}\|u(t)-w(t)\|
 \leq C_R\|\psi^{-1/2}R^c[w]\|_{L^2(Q_T)}.
 \label{eq:residual_bound}
\end{equation}
\end{corollary}
\begin{proof}
For two solutions with the same coefficients, data and kernel,
Theorem~\ref{thm:stability} applies with $f=0$ and $e(0)=0$.
For residual control, $e=u-w$ satisfies \eqref{eq:error_equation}
with $f=-R^c[w]$, zero initial error and zero boundary trace.
The same theorem gives \eqref{eq:residual_bound}.
\end{proof}

The independence from memory strength concerns the coefficient $C_R$.
The residual still depends on the memory law. Applying this estimate
to a grid calculation would require a reconstruction with the stated
regularity and control of its continuous residual, including
discretization defects.

With constant $\alpha=a$, constant $c,\nu$ and a monotone reaction, $b=0$ and $C_R=\sqrt{T^a/a}$. No inequality between $a$ and $\gamma(0)$ is required. For comparison, applying the elementary symmetric bound to the history
inner product and then discarding its remaining positive current-state
contribution leads to
\begin{equation}\label{eq:legacy_constant}
 C_{\rm old}=\exp\left\{\tfrac12\int_0^T\psi^{-1}
  [2b+1+\eta k_\rho(t)]\dd t\right\}.
\end{equation}
This elementary bound can diverge when $a\le\gamma(0)$. The running-maximum estimate avoids that restriction by keeping the present and past contributions together. The comparison is limited to the displayed elementary bound; it does not cover all existing nonlocal stability estimates.

\subsection{A sharp refinement using local coercivity}\label{sec:coercivity}
The preceding estimate discards the positive diffusion term. Retaining it improves the residual coefficient when the local operator is coercive. Suppose that the compared states satisfy, almost everywhere,
\begin{equation}\label{eq:local_coercivity}
 \langle c\cdot\nabla e-\nu\Delta e+\mu[\phi(u)-\phi(v)],e\rangle
 \ge\chi\|e\|^2,\qquad \chi>0.
\end{equation}
If $\|e\|\le C_P\|\nabla e\|$, the estimates above give the sufficient choice
$\chi=\nu_0/(2C_P^2)-b$ whenever it is positive. For spatially constant $\nu$, the full diffusion contribution is retained and one may use
$\chi=\nu_0/C_P^2-\tfrac12\|\nabla\!\cdot c\|_\infty-\kappa_\phi>0$.

\begin{theorem}[Coercivity-refined residual stability]\label{thm:coercive}
Let Assumption~\ref{ass:standing} hold. Suppose that
\eqref{eq:local_coercivity} holds with a constant $\chi>0$
and that $\psi^{-1/2}f\in L^2(Q_T)$. Then
\begin{equation}\label{eq:coercive_bound}
 \max_{0\le t\le T}\|e(t)\|
 \le\|e(0)\|+C_\chi(T)\|\psi^{-1/2}f\|_{L^2(Q_T)},\qquad
 C_\chi(T)=\sqrt{\frac{1-e^{-2\chi H(T)}}{2\chi}}.
\end{equation}
For fixed $\chi$ and clock $H$, the coefficient is independent of
memory strength and kernel mass. For $T>0$,
$0<C_\chi(T)<\min\{\sqrt{H(T)},(2\chi)^{-1/2}\}$.
For zero initial error, $C_\chi(T)$ is the optimal uniform coefficient
depending only on $\chi$ and $H(T)$ over the stated class of operators
and admissible memory terms, including $\eta=0$.
\end{theorem}
\begin{proof}
Use $y$ and $Y$ from Theorem~\ref{thm:stability}. Fix $0<\varepsilon<T$
and define, up to null sets,
\[
 A_\varepsilon=\{t\in(\varepsilon,T):Y'(t)>0\},\qquad
 h_\varepsilon(t)=\int_\varepsilon^t
       \boldsymbol1_{A_\varepsilon}(r)\psi(r)^{-1}\dd r.
\]
Almost everywhere on $A_\varepsilon$, $y=Y>0$ and $y'=Y'$.
The history inner product is nonnegative there.
Testing \eqref{eq:error_equation} and using
\eqref{eq:local_coercivity} gives
$Y'+\chi\psi^{-1}Y\le\psi^{-1}\|f\|$ on this set.
Since $Y'=0$ almost everywhere outside $A_\varepsilon$,
\[
 Y'+\chi h_\varepsilon'Y\le h_\varepsilon'\|f\|
 \quad\text{almost everywhere on }[\varepsilon,T].
\]
The integrating factor $e^{\chi h_\varepsilon}$ gives
\[
 Y(T)\le e^{-\chi h_\varepsilon(T)}Y(\varepsilon)
 +\int_\varepsilon^T
 e^{-\chi[h_\varepsilon(T)-h_\varepsilon(r)]}
 h_\varepsilon'(r)\|f(r)\|\dd r.
\]
Apply Cauchy--Schwarz with respect to $h_\varepsilon'(r)\dd r$.
The chain rule gives
\[
 \int_\varepsilon^T
 e^{-2\chi[h_\varepsilon(T)-h_\varepsilon(r)]}
 h_\varepsilon'(r)\dd r
 =\frac{1-e^{-2\chi h_\varepsilon(T)}}{2\chi},
\]
even when $h_\varepsilon$ has flat portions. Since
$0\le h_\varepsilon(T)\le H(T)-H(\varepsilon)$ and
$\int_\varepsilon^T h_\varepsilon'\|f\|^2\dd r
 \le\|\psi^{-1/2}f\|_{L^2(Q_T)}^2$, it follows that
\[
 Y(T)\le Y(\varepsilon)+C_\chi(T)\|\psi^{-1/2}f\|_{L^2(Q_T)}.
\]
Letting $\varepsilon\downarrow0$ proves \eqref{eq:coercive_bound}.
The coefficient bounds follow from
$C_\chi(T)^2=\int_0^{H(T)}e^{-2\chi s}\dd s$.

For sharpness, take $\Omega=(0,1)$, $c=0$, $\eta=\mu=0$ and
$\nu=\chi/\pi^2$. The Dirichlet Poincar\'e inequality gives
\[
 \langle-\nu v_{xx},v\rangle
 =\nu\|v_x\|^2\ge\nu\pi^2\|v\|^2=\chi\|v\|^2,
 \qquad v\in H^2(0,1)\cap H_0^1(0,1).
\]
Thus this operator satisfies \eqref{eq:local_coercivity}.
Let $X(x)=\sqrt2\sin(\pi x)$, so $\|X\|=1$ and
$-\nu X''=\chi X$. With $L=H(T)$, set
\[
 f(t)=e^{-\chi[L-H(t)]}X,\qquad
 e(t)=q(H(t))X,\qquad
 q(h)=e^{-\chi L}\frac{\sinh(\chi h)}{\chi}.
\]
The identity $q'+\chi q=e^{-\chi(L-h)}$ proves the error equation.
The compared states $e$ and $0$ have the required regularity and
zero initial and boundary data. Since $q'>0$,
\[
 \max_t\|e(t)\|=q(L)=C_\chi(T)^2,\qquad
 \|\psi^{-1/2}f\|_{L^2(Q_T)}^2
 =\int_0^L e^{-2\chi(L-h)}\dd h=C_\chi(T)^2.
\]
Hence equality holds, and no smaller uniform residual coefficient
is possible.
\end{proof}

The restriction to $A_\varepsilon$ is essential: a running maximum
can remain constant while the current error decreases.
Under \eqref{eq:local_coercivity}, $C_\chi(T)$ replaces $C_R$ in
Corollary~\ref{cor:error}. Its upper bound is uniform in $T$ when
the same positive $\chi$ applies; the forcing norm may still grow.
Sharpness over the full class does not imply optimality for each
fixed spatial operator or positive memory coupling.

\subsection{Residual control for the transport application}
\label{sec:app-residual-theory}
\begin{proposition}[Residual control for the porous-medium application]
\label{prop:app-residual-control}
Let $Q=(0,1)\times(0,1)$ and use the coefficients of
\eqref{eq:app-conformable}. Assume existence of an exact state $u$.
Let $u$ and a trial state $w$ have the solution regularity specified in
Assumption~\ref{ass:standing} and identical initial and boundary data.
With $a(t)=\kappa_d+\ell(H(t))$, define
\[
 R^c[w]=\psi w_t+(v-D')w_x-Dw_{xx}
       +a(t)w+\eta(t)\Cc_\rho[w]-q(x,H(t)).
\]
If $\psi^{-1/2}R^c[w]\in L^2(Q)$, then
\begin{equation}
 \max_{0\leq t\leq1}\|u(t)-w(t)\|_{L^2(0,1)}
 \leq\sqrt{H(1)}\,\|\psi^{-1/2}R^c[w]\|_{L^2(Q)}.
 \label{eq:app-residual-control}
\end{equation}
\end{proposition}
\begin{proof}
Set $e=u-w$. The conservative structure and zero boundary trace give
\begin{equation}
 \langle(v-D')e_x-De_{xx}+a(t)e,e\rangle
 =\int_0^1D|e_x|^2\,dx+a(t)\|e\|_2^2\geq0.
 \label{eq:app-energy}
\end{equation}
The $D'$ terms cancel before estimation. At an active maximum of the
spatial error norm, the memory contribution is nonnegative because
$\eta(t)\geq0$ and $\rho\geq0$. The argument of
Theorem~\ref{thm:stability} therefore gives
$Y'\leq\psi^{-1}\|R^c[w]\|_2$ almost everywhere away from zero.
Integrate from $\varepsilon$ to $1$, pass to zero using the assumed
regularity and $Y(0)=0$, and apply Cauchy--Schwarz with
$\int_0^1\psi^{-1}\,dt=H(1)$.
\end{proof}
This proof allows the prescribed nonnegative coupling $\eta(t)$ and
the time-dependent monotone local term $a(t)u$. Existence and state
regularity remain assumptions. The estimate concerns the continuous
spatial norm and integrated residual; it does not itself prove
convergence of either compartment time-stepping scheme.

\subsection{Analysis of the compartment schemes}
\label{sec:p3pair-analysis}
The following statements concern the finite-dimensional system defined
in Section~\ref{sec:method}. Assume $D>0$, $\kappa_d\geq0$, positive
rates, nonnegative capacities, and the mesh condition
\eqref{eq:p3pair-mesh}. They do not follow solely from the continuous
residual estimates.

\begin{proposition}[Solvability of the implicit steps]
\label{prop:p3pair-implicit-solvability}
For every $\tau>0$, the implicit mobile system
\eqref{eq:p3pair-implicit} has a unique solution. Its coefficient matrix
is a nonsingular M-matrix. Compartment recovery is unique.
\end{proposition}
\begin{proof}
The off-diagonal entries of $L_h$ are nonnegative. Applying the flux sum
to each coordinate vector gives column sums
$-2D(0)/h^2$ and $-2D(1)/h^2$ at the first and last columns,
respectively, and zero at interior columns. Thus the column sums are
nonpositive. Let $C_n=(a_0+\kappa_d+\chi_n)I-L_h$. Its off-diagonal
entries are nonpositive, and the diagonal-dominance margin of column
$j$ is
\[
 (C_n)_{jj}-\sum_{i\ne j}|(C_n)_{ij}|
 =a_0+\kappa_d+\chi_n-\sum_i(L_h)_{ij}>0.
\]
To see inverse positivity directly, write $C_n=D_n-B_n$, where $D_n$
is its positive diagonal and $B_n\geq0$ has zero diagonal. Strict column
dominance gives $\|B_nD_n^{-1}\|_1<1$. Hence
$C_n^{-1}=D_n^{-1}\sum_{m=0}^{\infty}(B_nD_n^{-1})^m\geq0$.
Finally, $a_0+\lambda_k>0$ in \eqref{eq:p3pair-recovery}.
\end{proof}
The backward-Euler startup is positive for nonnegative initial data and
source. For subsequent BDF2 steps, inverse positivity does not make the
history terms nonnegative: $B^n$ and $P_k^n$ contain negative multiples
of the older state. Unconditional BDF2 positivity is therefore not
asserted.

\begin{proposition}[Explicit positivity and weighted stability]
\label{prop:p3pair-explicit-stability}
If $0<\tau\leq\tau_{\rm pos}$, forward Euler preserves nonnegative
mobile and compartment states for nonnegative initial data and source.
For positive capacities define
\[
 \|Y\|_{1,w}
   =h\sum_i|U_i|+h\sum_{i,k}b_k|Z_{k,i}|,\qquad
 Y=(U,Z_1,\ldots,Z_{N_c}).
\]
For two explicit trajectories with error $E^n$ and mobile source
difference $F^n$,
\begin{equation}
 \|E^n\|_{1,w}
 \leq\|E^0\|_{1,w}
  +\tau h\sum_{m=0}^{n-1}\sum_i|F_i^m|.
 \label{eq:p3pair-explicit-stability}
\end{equation}
In particular, the unforced update is nonexpansive in this norm.
For the no-trapping case, omit the uncoupled immobile variables.
\end{proposition}
\begin{proof}
Write the semidiscrete system as $\dot Y=AY+(Q,0,\ldots,0)$.
Its off-diagonal entries are nonnegative. Its mobile diagonal entries
are $(L_h)_{ii}-\kappa_d-S_{N_c}$ and its compartment diagonal entries are
$-\lambda_k$. The restriction \eqref{eq:p3pair-positive-step} therefore
makes $P=I+\tau A$ entrywise nonnegative. Nonnegative data and source
remain nonnegative by induction.

Let $w$ contain $h$ for mobile entries and $hb_k$ for compartment $k$.
Exchange cancels in $w^TAY$. Its remaining coefficients are the
nonpositive transport column sums and mobile degradation, so
$w^TA\leq0$ componentwise. Consequently $w^TP\leq w^T$ and
\[
 \|PY\|_{1,w}=w^T|PY|\leq w^TP|Y|\leq w^T|Y|.
\]
Apply this inequality to
$E^{n+1}=PE^n+\tau(F^n,0,\ldots,0)$ and sum over the steps.
\end{proof}

\begin{proposition}[Discrete mass balance]
\label{prop:p3pair-balance}
With the synchronized cumulative updates of
Section~\ref{sec:p3pair-accounting}, both schemes satisfy
\begin{equation}
 M_m^n+M_i^n+M_{\rm cap}^n+M_{\rm deg}^n
 =M_m^0+M_i^0+I_{\rm inj}^n
 \label{eq:p3pair-balance}
\end{equation}
in exact arithmetic.
\end{proposition}
\begin{proof}
Multiply each mobile equation by $h$ and each compartment equation
by $hb_k$, then sum. The exchange contributions cancel pairwise.
The transport terms telescope to $-\Phi$, and the remaining rates
are $h\sum_iQ_i-\kappa_dM_m-\Phi$. Adding the cumulative equations
therefore gives zero time difference for the balance defect.
For forward Euler this preserves its initial value at every step.
For the implicit method, backward Euler gives the same defect at
levels zero and one; the BDF2 recurrence
$3d^n-4d^{n-1}+d^{n-2}=0$ then preserves it by induction.
The initial defect is zero.
\end{proof}
The balance identity measures algebraic conservation under synchronized
time integration. Replacing the discrete injection by its exact integral
also reveals the source-integration discrepancy.
The nxt subsection estimates the state error separately
\subsection{Combined spatial, compartment and temporal error}
\label{sec:p3ce-combined}
We now compare the fully discrete mobile state with cell averages of the
continuous-rate PDE solution. We write $M=N_c$ for the compartment count
in the error analysis. All times in this subsection are physical
times $\vartheta\in[0,\Theta]$. All temporal derivatives below
refer to this coordinate. Write
\begin{equation}
 Lu=-\partial_x(vu-Du_x),\qquad
 u_\vartheta=Lu-\kappa_du-\ell*u_\vartheta+f,
 \label{eq:p3ce-continuous}
\end{equation}
on $(0,1)$, with homogeneous Dirichlet values and initially empty mobile
and immobile phases. Here $*$ is causal convolution,
$\ell(\vartheta)=\int b(\lambda)\lambda e^{-\lambda\vartheta}\,d\lambda$,
with $b(\lambda)\ge0$, $\kappa_d\ge0$, and
\[
 B=\int b(\lambda)\,d\lambda<\infty,\qquad
 S=\int b(\lambda)\lambda\,d\lambda<\infty.
\]
The source $f$ denotes either the prescribed application source or the
manufactured source. The following result is conditional on the stated
continuum solution regularity; it does not establish that regularity
from the application data.

Let $h=1/N_x$, $N_x\ge2$, and let $\Pi_h$ denote exact cell averaging.
Use the conservative flux operator $L_h$ of
Section~\ref{sec:method}, and set
$\|V\|_{2,h}^2=h\sum_i|V_i|^2$. Assume $v$ is constant,
$D\in W^{1,\infty}(0,1)$, $0<d_*\le D\le d^*$, and the sufficient
all-face condition
\begin{equation}
 |v|h\le2D(jh),\qquad j=0,\ldots,N_x
 \label{eq:p3ce-mesh}
\end{equation}
For the application coefficient, $D(0)=D(1)=\max D$, so the implemented
interior-face restriction also implies the two endpoint conditions.
Put
\[
 X=\{w\in W^{3,\infty}(0,1):w(0)=w(1)=0\},
 \qquad u\in W^{1,1}(0,\Theta;X),\quad u(\cdot,0)=0.
\]
Assume $Q_h=\Pi_hf\in W^{2,1}(0,\Theta;\mathbb R^{N_x})$ and $Q_h(0)=0$.
Positive capacities $b_j$ and rates $\lambda_j$ define
$\ell_M(\vartheta)=\sum_{j=1}^M b_j\lambda_j e^{-\lambda_j\vartheta}$.
We distinguish $U_h$, the continuous-time solution with $\ell$;
$U_{h,M}$, the continuous-time solution with $\ell_M$; and
$U_{h,M}^n$, the numerical mobile state at $\vartheta_n=n\tau$,
$\tau=\Theta/N_\vartheta$. Each uses the same $L_h$ and $Q_h$.

\paragraph*{Spatial projection and propagation.}\enspace 
Summation by parts for the actual boundary fluxes gives
\begin{align}
 -\langle L_hV,V\rangle_h
 &=\sum_{j=1}^{N_x-1}\frac{D(jh)}h(V_{j+1}-V_j)^2\notag\\
 &\quad+\left(\frac{2D(0)}h+\frac v2\right)V_1^2
       +\left(\frac{2D(1)}h-\frac v2\right)V_{N_x}^2.
 \label{eq:p3ce-energy}
\end{align}
If $|V|_{1,h}^2=h^{-1}[V_1^2+V_{N_x}^2+
\sum_j(V_{j+1}-V_j)^2]$, then \eqref{eq:p3ce-mesh} and telescoping imply
$-\langle L_hV,V\rangle_h\ge d_*|V|_{1,h}^2$ and
$\|V\|_{2,h}\le|V|_{1,h}$. Thus $L_h$ is invertible and the elliptic
projection $R_hw$ defined by $L_hR_hw=\Pi_hLw$ is well defined.

\begin{lemma}[Cell-average projection estimate]
\label{lem:p3ce-projection}
For $w\in X$,
\begin{equation}
 \|R_hw-\Pi_hw\|_{2,h}\le h^2A(w),\qquad
 A(w)=c_2\|w_{xx}\|_\infty+c_3\|w_{xxx}\|_\infty,
 \label{eq:p3ce-projection}
\end{equation}
where
\[
 c_2=\frac16+\frac{|v|}{3d_*}+\frac{2d^*}{3d_*},
 \qquad c_3=\frac{d^*}{12d_*}.
\]
\end{lemma}
\begin{proof}
Let $F_j$ be the implemented face flux and set
$\delta_j=F_j(\Pi_hw)-(vw-Dw_x)(jh)$. Taylor's integral formula and
the symmetric cell-average kernels give
\begin{align*}
 |\delta_j|&\le h^2\left(\frac{|v|}{6}\|w_{xx}\|_\infty+
 \frac{d^*}{12}\|w_{xxx}\|_\infty\right)
 \quad(1\le j<N_x),\\
 |\delta_0|&\le\frac{D(0)h}{3}\|w_{xx}\|_\infty,
 \qquad |\delta_{N_x}|\le\frac{D(1)h}{3}\|w_{xx}\|_\infty.
\end{align*}
In particular, the boundary face defect is only $O(h)$; a uniform
$O(h^2)$ cell-residual argument is not available. with
$\eta=R_hw-\Pi_hw$ and $(\operatorname{div}_h\delta)_i
=(\delta_i-\delta_{i-1})/h$, exact cell balance gives
$L_h\eta=\operatorname{div}_h\delta$. Choose $W$ affine through the cell
centres with $W_1=h\delta_0/(2D(0))$ and
$W_{N_x}=-h\delta_{N_x}/(2D(1))$. Then
\[
 \|W\|_\infty\le h^2\|w_{xx}\|_\infty/6,
 \qquad |\operatorname{slope}W|
 \le h^2\|w_{xx}\|_\infty/[3(1-h)].
\]
The vector $r=\eta-W$ satisfies
$L_hr=\operatorname{div}_h d$ with $d=\delta+F(W)$ and $d_0=d_{N_x}=0$.
Its interior face values obey
\[
 |d_j|\le h^2\left[
 \left(\frac{|v|}{3}+\frac{d^*}{3(1-h)}\right)\|w_{xx}\|_\infty
 +\frac{d^*}{12}\|w_{xxx}\|_\infty\right].
\]
Summation by parts, coercivity and Cauchy--Schwarz yield
$\|r\|_{2,h}\le\max_j|d_j|/d_*$. Add the bound for $W$ and use
$(1-h)^{-1}\le2$.
\end{proof}

For a zero-initial-data spatial system with continuous memory and mobile
forcing $g$, the following contractive forcing estimate holds:
\begin{equation}
 \max_{\vartheta\le\Theta}\|E(\vartheta)\|_{2,h}
 \le\int_0^\Theta\|g(s)\|_{2,h}\,ds.
 \label{eq:p3ce-forcing}
\end{equation}
Indeed $K=-\ell'\ge0$ belongs to $L^1(0,\Theta)$ and
$\int_0^\vartheta K=S-\ell(\vartheta)$. The error equation is
$E'=(L_h-\kappa_dI-SI)E+K*E+g$. Dissipativity gives the norm inequality
$m'\le-(\kappa_d+S)m+K*m+\|g\|_{2,h}$.
At an active running maximum of $m=\|E\|_{2,h}$, the first two terms
are at most $-(\kappa_d+\ell(\vartheta))m\le0$.
Integrating the running-maximum inequality proves \eqref{eq:p3ce-forcing}.

Set $\eta=R_hu-\Pi_hu$. The residual of $R_hu$ in the spatial system is
$\eta'+\kappa_d\eta+\ell*\eta'$. Using
$\|\ell\|_{L^1(0,\Theta)}\le B$ in \eqref{eq:p3ce-forcing}, followed by
Lemma~\ref{lem:p3ce-projection}, gives
\begin{align}
 \max_{\vartheta\le\Theta}\|U_h-\Pi_hu\|_{2,h}
 &\le h^2C_{\rm sp}(u,\Theta),\label{eq:p3ce-space}\\
 C_{\rm sp}(u,\Theta)
 &=\sup_{\vartheta\le\Theta}A(u(\vartheta))
 +(1+B)\int_0^\Theta A(u_\vartheta(s))\,ds
 +\kappa_d\int_0^\Theta A(u(s))\,ds.
 \label{eq:p3ce-space-constant}
\end{align}
Neither an inverse power of $h$ nor a largest exchange rate appears in
this constant; the solution regularity is essential.

\begin{remark}
	The term $|v|/(3d_*)$ in $c_2$ contributes to the conservatism
	of the analytical projection bound. For a prescribed profile
	$\varphi$, the projection error can be evaluated directly
	from the tridiagonal system
	$L_hR_h\varphi=\Pi_hL\varphi$.
	This allows the spatial contribution to be assessed
	without replacing the projection error by its general bound.
\end{remark}

\paragraph*{Compartment approximation on the same mesh.}\enspace 
Define
\begin{align}
 D_M(\vartheta)&=\int_0^\vartheta[\ell(s)-\ell_M(s)]\,ds,
 &\varepsilon_M^{(1)}(\Theta)&=\int_0^\Theta|D_M(s)|\,ds,
 \label{eq:p3ce-memory-defect}\\
 C_{2,h}(\Theta)&=\Theta\|Q_h'(0)\|_{2,h}
 +\int_0^\Theta(\Theta-s)\|Q_h''(s)\|_{2,h}\,ds.
 \label{eq:p3ce-source-factor}
\end{align}
For the finite coupled state $Y=(U,Z_1,\ldots,Z_M)$, write
$Y'=A_MY+\mathcal F_h$ with $\mathcal F_h=(Q_h,0,\ldots,0)$ and use
\[
 \|Y\|_{\mathcal H_M}^2=\|U\|_{2,h}^2+
 \sum_jb_j\|Z_j\|_{2,h}^2.
\]
Exchange cancellation yields
\begin{equation}
 \langle A_MY,Y\rangle_{\mathcal H_M}
 =\langle L_hU,U\rangle_h-\kappa_d\|U\|_{2,h}^2
 -\sum_jb_j\lambda_j\|U-Z_j\|_{2,h}^2\le0.
 \label{eq:p3ce-coupled-dissipation}
\end{equation}
Hence $e^{\vartheta A_M}$ is contractive in $\mathcal H_M$.
Since $Y(0)=\mathcal F_h(0)=0$,
\[
 Y''(\vartheta)=e^{\vartheta A_M}\mathcal F_h'(0)
 +\int_0^\vartheta e^{(\vartheta-s)A_M}\mathcal F_h''(s)\,ds,
 \qquad
 \int_0^\Theta\|U_{h,M}''\|_{2,h}\,ds\le C_{2,h}(\Theta).
\]
The continuum-memory residual of $U_{h,M}$ is
$(\ell-\ell_M)*U_{h,M}'=D_M*U_{h,M}''$, where the initial term vanishes
because $U_{h,M}'(0)=Q_h(0)=0$. Young's inequality and
\eqref{eq:p3ce-forcing} therefore imply
\begin{equation}
 \max_{\vartheta\le\Theta}\|U_{h,M}-U_h\|_{2,h}
 \le C_{2,h}(\Theta)\varepsilon_M^{(1)}(\Theta)
 \label{eq:p3ce-memory}
\end{equation}
The prefactor is independent of $M$ and the largest rate. It is uniform
in $h$ when the displayed source norms are uniformly bounded.

For a computable estimate of $\varepsilon_M^{(1)}$, take a partition
$0=s_0<\cdots<s_J=\Theta$ and the linear interpolant $I_D$ of $D_M$.
Since $D_M''=-K+K_M$ is integrable, the interpolation Green kernel gives
\begin{equation}
 \varepsilon_M^{(1)}(\Theta)
 \le\int_0^\Theta|I_D(s)|\,ds+\sum_{j=0}^{J-1}\mathcal R_j,
 \qquad
 \mathcal R_j=\frac{\Delta_j^2}{8}
 [\ell(s_j)-\ell(s_{j+1})+\ell_M(s_j)-\ell_M(s_{j+1})],
 \label{eq:p3ce-kernel-enclosure}
\end{equation}
where $\Delta_j=s_{j+1}-s_j$. The integral of $|I_D|$ is evaluated
exactly on each segment, splitting at a zero when necessary.
If $s_j>0$, the positive exponential-mixture kernels $K$ and $K_M$
are decreasing, so one may use the smaller remainder
$\min\{\mathcal R_j,\Delta_j^3 H_j/12\}$, with
\[
 H_j=\max\{|K(s_j)-K_M(s_{j+1})|,
             |K(s_{j+1})-K_M(s_j)|\}.
\]
The first remainder remains valid at the integrably singular endpoint.
These are exact-arithmetic inequalities. A certified floating-point
enclosure would additionally require bounds on kernel and node-value
evaluation errors.

\paragraph*{A source-variation bound for the single-cell check.}\enspace
The scalar experiment below admits a sharper source factor than the
general transport estimate. Omit $L_h$ and let $u$ and $u_M$ solve
$u'+\kappa_du+\ell*u'=q$ and its positive $M$-compartment approximation,
with both phases initially empty. Assume $q\ge0$,
$q\in W^{1,1}(0,\Theta_c)$, $q(0)=0$, $\kappa_d\ge0$ and $S_M\le S$.
Then
\begin{equation}
 \max_{0\le\vartheta\le\Theta_c}|u(\vartheta)-u_M(\vartheta)|
 \le\left[\operatorname{TV}_{[0,\Theta_c]}(q)
       +2(\kappa_d+S)I(\Theta_c)\right]
       \varepsilon_M^{(1)}(\Theta_c),\qquad
 I(\Theta_c)=\int_0^{\Theta_c}q(s)\,ds.
 \label{eq:p3checks-single-cell-bound}
\end{equation}
Indeed, the finite exchange matrix is symmetric negative semidefinite
after weighting the compartment coordinates by $\sqrt{b_j}$.
Its mobile impulse response has the representation
$r_M(\vartheta)=\sum_k a_k e^{-\mu_k\vartheta}$ with
$a_k\ge0$, $\mu_k\ge0$ and $\sum_k a_k=1$.
Consequently $r_M'\le0$, $r_M''\ge0$ and
$r_M'(0)=-(\kappa_d+S_M)$. From $u_M=r_M*q$,
\[
 u_M''=q'+r_M'(0)q+r_M''*q,\qquad
 \|u_M''\|_{L^1(0,\Theta_c)}
 \le\operatorname{TV}(q)+2(\kappa_d+S_M)I(\Theta_c).
\]
The continuum residual is $D_M*u_M''$, since $u_M'(0)=q(0)=0$.
The contractive forcing estimate and Young's inequality give
Eq.~\eqref{eq:p3checks-single-cell-bound}. This scalar impulse-response
argument is specific to the single-cell check; its source factor is not
substituted for $C_{2,h}$ in the transport estimates.

\paragraph*{The actual time discretizations.}\enspace 
Put $S_M=\sum_jb_j\lambda_j$ and
\begin{equation}
 \tau_{\rm FE}=\left[
 \max\left\{\max_i[-(L_h)_{ii}+\kappa_d+S_M],\max_j\lambda_j\right\}
 \right]^{-1}.
 \label{eq:p3ce-fe-step}
\end{equation}
Under \eqref{eq:p3ce-mesh}, $L_h$ is Metzler with nonpositive column
sums and $L_h\mathbf1\le0$. For $\tau\le\tau_{\rm FE}$,
$P=I+\tau A_M\ge0$, $P\mathbf1\le\mathbf1$ and $P^Tw\le w$ for
$w=(h\mathbf1,hb_1\mathbf1,\ldots,hb_M\mathbf1)$.
Weighted Cauchy--Schwarz gives
\[
 \sum_i w_i|(PV)_i|^2
 \le\sum_i w_i\Bigl(\sum_jP_{ij}\Bigr)\sum_jP_{ij}|V_j|^2
 \le\sum_j(P^Tw)_j|V_j|^2
 \le\sum_jw_j|V_j|^2.
\]
Thus $\|PV\|_{\mathcal H_M}\le\|V\|_{\mathcal H_M}$.
Here the indices run over all mobile and compartment entries.
The forward-Euler local defect is
$\int_{\vartheta_n}^{\vartheta_{n+1}}
(\vartheta_{n+1}-s)Y''(s)\,ds$.
Combining contraction with the preceding representation of $Y''$ gives
\begin{equation}
 \max_n\|Y(\vartheta_n)-Y^n\|_{\mathcal H_M}
 \le\mathcal T_{\rm FE}:=\frac{\tau}{2}C_{2,h}(\Theta).
 \label{eq:p3ce-fe-time}
\end{equation}
for completeness, after interchanging integrals the coefficient of
$\|\mathcal F_h''(r)\|$ for $r\in[\vartheta_j,\vartheta_{j+1}]$ is
$[(\vartheta_{j+1}-r)^2+(N_\vartheta-j-1)\tau^2]/2
\le\tau(\Theta-r)/2$; the initial derivative contributes
$\Theta\tau\|\mathcal F_h'(0)\|/2$.

For backward-Euler-started BDF2, assume
$Y\in W^{3,1}(0,\Theta;\mathcal H_M)$. No step restriction is needed
for the following error estimate:
\begin{equation}
 \max_n\|Y(\vartheta_n)-Y^n\|_{\mathcal H_M}
 \le\mathcal T_{\rm BDF2}:=
 \sqrt5\int_0^\tau s\|Y''(s)\|_{\mathcal H_M}\,ds
 +\frac{\sqrt5}{4}\tau^2
       \int_0^\Theta\|Y'''(s)\|_{\mathcal H_M}\,ds.
 \label{eq:p3ce-bdf-time}
\end{equation}
The standard $G$-stability argument~\cite{Dahlquist1978} uses
$\|(a,b)\|_G^2=\|a\|_{\mathcal H_M}^2+
\|2a-b\|_{\mathcal H_M}^2$ and the identity
\[
 \left\langle\tfrac32a-2b+\tfrac12c,a\right\rangle_{\mathcal H_M}
 =\tfrac14\{\|(a,b)\|_G^2-\|(b,c)\|_G^2
                    +\|a-2b+c\|_{\mathcal H_M}^2\}.
\]
Together with \eqref{eq:p3ce-coupled-dissipation}, it makes the
homogeneous two-level map contractive. Since
$\|(\tfrac32I-\tau A_M)^{-1}\|\le2/3$, a derivative defect $d_n$
adds at most $(2\sqrt5/3)\tau\|d_n\|$ to the $G$ norm.
The backward-Euler resolvent is contractive, and its first-step error
is bounded by $\int_0^\tau s\|Y''(s)\|\,ds$.
The BDF2 derivative defect has the Peano representation
\[
 d_n=-\int_0^{2\tau}w_\tau(s)Y'''(\vartheta_n-s)\,ds,
 \quad
 w_\tau(s)=
 \begin{cases}s-3s^2/(4\tau),&0\le s\le\tau,\\
 (2\tau-s)^2/(4\tau),&\tau\le s\le2\tau.
 \end{cases}
\]
These nonnegative kernels overlap at most twice, with pointwise sum
at most $3\tau/8$. Thus
$\tau\sum_{n\ge2}\|d_n\|\le(3/8)\tau^2\int_0^\Theta\|Y'''\|$.
The initial pair $(e^1,0)$ has $G$ norm $\sqrt5\|e^1\|$,
which proves \eqref{eq:p3ce-bdf-time}. This argument establishes
stability and error control, not unconditional BDF2 positivity.

\begin{theorem}[Combined cell-average mobile-state estimate]
\label{thm:p3ce-combined}
Consider the initially empty problem \eqref{eq:p3ce-continuous}.
Assume the coefficient and solution regularity stated at the start
of this subsection, including \eqref{eq:p3ce-mesh}.
Let all compartment capacities and rates be positive.
Use exact source cell averages
$Q_h\in W^{2,1}(0,\Theta;\mathbb R^{N_x})$, with $Q_h(0)=0$.
Forward Euler with $\tau\le\tau_{\rm FE}$ then satisfies
\eqref{eq:p3ce-combined-bound}, with
$\mathcal T_{\rm method}=\mathcal T_{\rm FE}$.
Backward-Euler-started BDF2 satisfies the same estimate with
$\mathcal T_{\rm method}=\mathcal T_{\rm BDF2}$ if
$Y\in W^{3,1}(0,\Theta;\mathcal H_M)$:
\begin{equation}
 \max_n\|U_{h,M}^n-\Pi_hu(\vartheta_n)\|_{2,h}
 \le h^2C_{\rm sp}(u,\Theta)
     +C_{2,h}(\Theta)\varepsilon_M^{(1)}(\Theta)
     +\mathcal T_{\rm method}.
 \label{eq:p3ce-combined-bound}
\end{equation}
\end{theorem}
\begin{proof}
Insert $U_{h,M}(\vartheta_n)$ and $U_h(\vartheta_n)$ between the
fully discrete state and $\Pi_hu(\vartheta_n)$. Apply the triangle
inequality, \eqref{eq:p3ce-space}, \eqref{eq:p3ce-memory} and the
appropriate temporal estimate. Projection of the coupled state onto
its mobile component has norm one in $\mathcal H_M$.
\end{proof}
The spatial and memory prefactors contain no inverse mesh powers or
largest exchange rate. The Euler prefactor has the same property;
its admissible step still depends on stiffness. The BDF2 stability
constants are universal, but uniform temporal accuracy additionally
requires uniform bounds on the displayed derivatives of $Y$.
Those derivative bounds must be checked for each source; they do not
follow merely from dissipativity.

\paragraph*{A uniform specialization for the manufactured benchmark.}\enspace 
For $u_\star=\vartheta^2\varphi$, with
$\varphi(x)=\sin(\pi x)+x(1-x)$, define
\[
 \Phi_h=\Pi_h\varphi,\quad \Psi_h=\Pi_hL\varphi,\quad
 d_h=L_h\Phi_h-\Psi_h,\quad
 q_\ell(\vartheta)=2\int_0^\vartheta(\vartheta-s)\ell(s)\,ds.
\]
Define $q_{\ell_M}$ by the same integral with $\ell_M$ in place of $\ell$.
The fixed continuous-kernel manufactured source is
$Q_h=[2\vartheta+\kappa_d\vartheta^2+q_\ell(\vartheta)]\Phi_h
-\vartheta^2\Psi_h$. In particular,
$Q_h'(0)=2\Phi_h$; the backward-Euler startup error does not vanish.
Set
\begin{align}
 \overline C_{2,h}
 &=2\Theta\|\Phi_h\|_{2,h}
 +\Theta^2\|\kappa_d\Phi_h-\Psi_h\|_{2,h}
 +q_\ell(\Theta)\|\Phi_h\|_{2,h},\label{eq:p3ce-manufactured-c2}\\
 C_{Q,h}&=2\|\kappa_d\Phi_h-\Psi_h\|_{2,h}
 +2S\|\Phi_h\|_{2,h},\label{eq:p3ce-manufactured-cq}\\
 \overline C_{3,h}
 &=2\Theta\|d_h\|_{2,h}
    +[\sqrt{2\Theta S}+4\Theta S]\|\Phi_h\|_{2,h}.
 \label{eq:p3ce-manufactured-c3}
\end{align}
If $S_M\le S$, then valid terms in \eqref{eq:p3ce-combined-bound} are
\begin{align}
 C_{2,h}&\le\overline C_{2,h},\qquad
 C_{\rm sp}(u_\star,\Theta)
 =A(\varphi)[(2+B)\Theta^2+\kappa_d\Theta^3/3],\notag\\
 \mathcal T_{\rm FE}&\le\tfrac12\tau\overline C_{2,h},\qquad
 \mathcal T_{\rm BDF2}
 \le\sqrt5\tau^2\left[
 \|\Phi_h\|_{2,h}+\frac{\tau C_{Q,h}}3
 +\frac{\overline C_{3,h}}4\right].
 \label{eq:p3ce-manufactured-time}
\end{align}
Here $A(\varphi)=c_2(\pi^2+2)+c_3\pi^3$.

To verify the temporal constants without a largest-rate factor, let
$Y_{\star,M}$ have mobile component $\vartheta^2\Phi_h$ and compartment
components
\[
 Z_{\star,j}(\vartheta)=
 \left[\vartheta^2-\frac{2\vartheta}{\lambda_j}
       +\frac{2(1-e^{-\lambda_j\vartheta})}{\lambda_j^2}\right]\Phi_h.
\]
The difference $W=Y-Y_{\star,M}$ has mobile forcing
$r=\vartheta^2d_h+(q_\ell-q_{\ell_M})\Phi_h$ and zero initial state.
It satisfies $r(0)=r'(0)=0$,
$r''(0)=2[d_h+(S-S_M)\Phi_h]$ and
$r'''=-2(K-K_M)\Phi_h\in L^1$, where $K_M=-\ell_M'$.
Contractivity and Fubini give
\[
 \int_0^\Theta\|W'''\|_{\mathcal H_M}\,ds
 \le2\Theta\|d_h+(S-S_M)\Phi_h\|_{2,h}
 +2\|\Phi_h\|_{2,h}\int_0^\Theta(\Theta-s)|K-K_M|(s)\,ds.
\]
Since $Z_{\star,j}'''=2\lambda_je^{-\lambda_j\vartheta}\Phi_h$,
Cauchy--Schwarz in time yields
\[
 \int_0^\Theta\|Y_{\star,M}'''\|_{\mathcal H_M}\,ds
 \le\sqrt{2\Theta S_M}\|\Phi_h\|_{2,h}.
\]
Use $\int_0^\Theta K\le S$, $\int_0^\Theta K_M\le S_M$ and
$S_M\le S$ to obtain
$\int_0^\Theta\|Y'''\|_{\mathcal H_M}\le\overline C_{3,h}$.
Furthermore $\|Q_h''(s)\|_{2,h}\le C_{Q,h}$, so the representation of
$Y''$ gives
$\int_0^\tau s\|Y''(s)\|_{\mathcal H_M}\,ds
\le\tau^2\|\Phi_h\|_{2,h}+\tau^3C_{Q,h}/3$.
This proves \eqref{eq:p3ce-manufactured-time}, including its nonzero
startup contribution.

These prefactors are uniform in both mesh and compartment resolution.
Indeed cell averaging is an $L^2$ contraction, and the preceding face
defects imply
\[
 \|d_h\|_{2,h}\le2h a_\varphi+
       \frac{d^*\sqrt{2h}}3\|\varphi''\|_\infty,
 \qquad
 a_\varphi=\frac{|v|}{6}\|\varphi''\|_\infty
 +\frac{d^*}{12}\|\varphi'''\|_\infty.
\]
The Gauss--Jacobi compartment rule used here satisfies $S_M\le S$.
Under $r=\lambda_0/\lambda$, $S$ integrates $\lambda_0/r$ against the
positive beta capacity measure. Gaussian quadrature underestimates
$1/(r+\epsilon)$ because its derivative of order $2M$ is positive.
Letting $\epsilon\downarrow0$ gives the assertion, since $S<\infty$.
For another compartment rule, $S_M\le S$ must be checked or replaced
by a proved uniform first-moment bound.

\paragraph*{Application source and observable scope.}\enspace 
For the prescribed delayed pulse $Q_h=P g_h$, with
$P(\vartheta)=10\sin^2[\pi(\vartheta-0.05)/0.2]$ on $[0.05,0.25]$
and zero otherwise, $\Theta\ge0.25$ gives
\[
 C_{2,h}(\Theta)=200\pi(\Theta-0.15)\|g_h\|_{2,h}.
\]
Although $P''$ jumps at the pulse endpoints, its distributional
derivative is a finite signed measure of total variation $3000\pi^2$,
symmetric about $0.15$. Because $P(0)=P'(0)=P''(0+)=0$, the Stieltjes
variation-of-constants formula gives
\[
 \int_0^\Theta\|Y'''\|_{\mathcal H_M}\,ds
 \le3000\pi^2(\Theta-0.15)\|g_h\|_{2,h}.
\]
If $\tau\le0.05$, the first exact and backward-Euler states are both
zero. The startup term in \eqref{eq:p3ce-bdf-time} then vanishes for
this pulse only. Jensen's inequality gives
$\|g_h\|_{2,h}\le\|g\|_{L^2(0,1)}$, so these application temporal and
memory prefactors are mesh- and rate-uniform. The application spatial
term still requires the stated continuum regularity and is not assigned
a computed value from the refinement differences.

The theorem controls absolute mobile-state error at the time nodes.
It also controls mobile total mass, since
\[
 \left|h\sum_i U_{h,M,i}^n-\int_0^1u(x,\vartheta_n)\,dx\right|
 \le\|U_{h,M}^n-\Pi_hu(\vartheta_n)\|_{2,h}.
\]
For the continuous piecewise-linear reconstruction $I_hV$ through
$(0,0),(x_i,V_i),(1,0)$,
\[
 \|I_hV-u\|_{L^2(0,1)}
 \le\|V-\Pi_hu\|_{2,h}+\frac{h^2}{6}\|u_{xx}\|_\infty.
\]
Here $\|I_hV\|_{L^2}\le\|V\|_{2,h}$, and the additional term follows
from centre interpolation and exact cell averaging. A piecewise-constant
reconstruction generally has an additional first-order approximation
error. No mesh-uniform point-detector, captured boundary-flux or
continuum immobile-mass estimate is asserted by
\eqref{eq:p3ce-combined-bound}. The corresponding application-output
refinement checks therefore retain their separate role.

\section{Numerical verification and reactive transport application}
\label{sec:results}
The experiments examine implementation, memory approximation, time
stepping and spatial resolution separately.
The manufactured test uses absolute cell-average errors.
The coupled-system and application comparisons use the relative norms
specified below. Each study identifies its reference state and
observation times before reporting results.

The refinement studies and application results use a completed
Windows run with Python 3.10.20, NumPy 2.2.6 and SciPy 1.15.3.
The inverse-Laplace calculations use mpmath 1.3.0.
Auxiliary reference checks retain separate provenance records.
The main computations use double precision and a numerical-library
CPU thread limit of one. Reference inversions use arbitrary precision.
The accompanying package records configurations, source hashes,
full-precision results and computed arrays.
Errors and orders are rounded only for presentation.
CPU times describe one run; they do not establish a general speed advantage.
Analytical bounds are evaluated in ordinary floating point.
No interval-arithmetic certification is claimed.

\subsection{Independent implementation and temporal checks}
\label{sec:p3pair-verification}
On a system with $100$ cells and $4$ compartments, we independently
assemble the full coupled matrix from direct face-flux differences.
Its transport action differs from the implementation by
$5.68\times10^{-14}$ on a fixed random vector with seed $23$.
One-step checks compare the explicit update and the eliminated implicit
update against the full coupled equations, including backward-Euler
startup and a BDF2 step. The largest state difference is
$4.44\times10^{-16}$. An explicit step above the implemented bound is
rejected. Additional nonzero-initial-state tests check mass accounting.

For an independent temporal comparison, the full coupled system is
integrated by adaptive BDF with relative tolerance $2\times10^{-11}$,
absolute tolerance $2\times10^{-13}$ and maximum step $5\times10^{-4}$.
This is an approximate semidiscrete reference, not an exact PDE solution.
We compare mobile cell values at $257$ common physical times. The norm
uses the spatial sum $h\sum_i|\cdot|^2$ and trapezoidal integration in
physical time. Table~\ref{tab:p3pair-temporal} reports the relative
discrepancy and the observed order
$p=\log_2(E_{\tau}/E_{\tau/2})$.

\begin{table}[htbp]
\centering\small
\caption{Temporal verification against the adaptive coupled-system
reference at fixed $N_x=100$ and $N_c=4$. Observed orders refer to
successive time-step halvings.}
\label{tab:p3pair-temporal}
\begin{tabular}{lrrr}
\toprule
Method & $N_\vartheta$ & Relative discrepancy & Observed order\\
\midrule
Implicit & 2048 & $1.2016\times10^{-4}$ & -- \\
Implicit & 4096 & $3.0042\times10^{-5}$ & 1.9999 \\
Implicit & 8192 & $7.5106\times10^{-6}$ & 2.0000 \\
\midrule
Explicit & 2048 & $1.3062\times10^{-2}$ & -- \\
Explicit & 4096 & $6.4578\times10^{-3}$ & 1.0163 \\
Explicit & 8192 & $3.2110 10^{-3}$ & 1.0080 \\
\bottomrule
\end{tabular}
\end{table}

The observed temporal orders are close to two for backward-Euler/BDF2
and to one for forward Euler. This fixed-system test does not estimate
spatial or compartment error.

\subsection{Separate memory and time--memory bound checks}
\label{sec:p3checks-separate-bounds}
The following two checks use the prescribed pulse and isolate different
parts of the analysis. The scalar experiment tests only the exchange-rate
approximation over a long observation horizon. The fixed-mesh experiment
then combines compartment and temporal errors for spatial transport.
Their norms and reference states are specified separately.
\subsubsection{Continuous-time single-cell memory bound}
\label{sec:p3checks-single-cell}
We first isolate the finite-spectrum approximation in a single cell,
without transport or time stepping. The mobile and immobile phases are
initially empty, with
$(\gamma,\lambda_0,\kappa_m,\kappa_d)=(0.5,0.5,0.8,0.1)$ and
$q(\vartheta)=10\sin^2[\pi(\vartheta-0.05)/0.2]$ on
$[0.05,0.25]$, extended by zero outside that interval.
The horizon is $\Theta_c=20$, which differs from the transport horizon
$\vartheta_f=1.0515132504186746$. This longer interval includes delayed
release after the pulse. Here $I(20)=1$ and
$\operatorname{TV}(q)=20$, so the source factor in
Eq.~\eqref{eq:p3checks-single-cell-bound} is
\begin{equation}
 20+2(\kappa_d+S)\approx22.4627416998
 \label{eq:p3checks-single-cell-factor}
\end{equation}
The primitive-kernel norm uses the interpolation remainder in
Eq.~\eqref{eq:p3ce-kernel-enclosure}.

Each finite-compartment state is evaluated by symmetric eigendecomposition
and analytical convolution with the pulse. A $128$-compartment trajectory
on a grid with $15881$ distinct times locates candidate error maxima;
the grid spacing is $5\times10^{-5}$ on $[0,0.6]$ and $0.005$ on
$[0.6,20]$. At each candidate maximum, direct numerical inverse Laplace
transformation supplies the continuous-spectrum value. Precision
refinement and independent de Hoog and Talbot inversions check the
reference at additional prescribed times. Thus the reported errors are
sampled discrepancies with continuum checks, not certified maxima over
the whole interval.

\begin{table}[htbp]
\centering\small
\caption{Single-cell memory approximation on $[0,20]$.
All errors are absolute mobile-concentration differences.
The upper estimate uses Eq.~\eqref{eq:p3checks-single-cell-bound};
the state discrepancy is sampled as described in the text.}
\label{tab:p3checks-single-cell}
\begin{tabular}{rrrr}
\toprule
$M$ & Sampled state error & Evaluated upper estimate & Estimate/error\\
\midrule
8 & $1.16380652\times10^{-4}$ & $1.09512544\times10^{-3}$ & 9.4\\
16 & $3.78947114\times10^{-6}$ & $4.20651405\times10^{-5}$ & 11.1\\
32 & $1.61688180\times10^{-8}$ & $1.44790807\times10^{-6}$ & 89.5\\
64 & $1.16621157\times10^{-11}$ & $4.71323688\times10^{-8}$ & 4041.5\\
\bottomrule
\end{tabular}
\end{table}

\begin{figure}[htbp]
\centering
\includegraphics[width=0.78\linewidth]{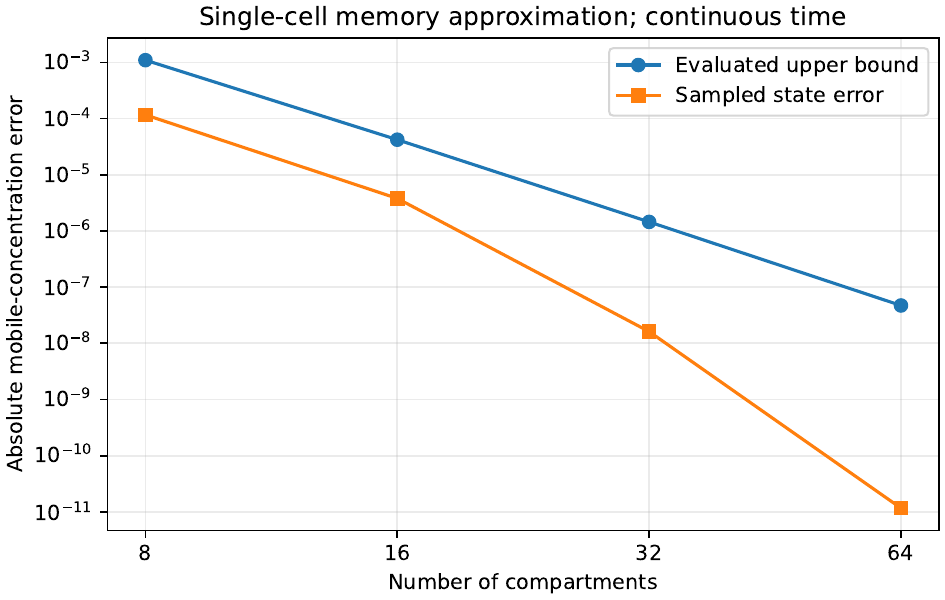}
\caption{Single-cell continuous-time memory check on $[0,20]$.
The evaluated analytical upper estimate and independently checked
sampled mobile-state discrepancy decrease with the compartment count.
Neither curve includes spatial discretization or time-stepping error.}
\label{fig:p3checks-memory-bound}
\end{figure}

Figure~\ref{fig:p3checks-memory-bound} shows that the estimate is
conservative, increasingly so at the larger compartment counts.
Among the tested values, its evaluated upper estimate selects
$M=16,16,32,64$ for absolute tolerances
$10^{-3},10^{-4},10^{-5},10^{-6}$, respectively.
These choices apply to this initially empty single-cell pulse.
The spatial transport problem requires its own compartment check.

\subsubsection{Pulse-source time and memory refinement on a fixed mesh}
\label{sec:p3checks-time-memory}
We next retain the spatial transport operator on $N_x=100$ cells and
use the application pulse source, with initially empty mobile and
immobile phases. This experiment isolates time and compartment errors
on a common mesh; it does not include spatial PDE error.
Let $U_h$ denote the continuous-time, continuous-rate state and
$U_{h,M}$ its continuous-time $M$-compartment approximation on that mesh.
For the numerical state $U_{h,M}^n$, set
\begin{align}
 E_{\rm time}&=\max_{\vartheta_n\in\mathcal T}
       \|U_{h,M}^n-U_{h,M}(\vartheta_n)\|_{2,h},\nonumber\\
 E_{\rm tm}&=\max_{\vartheta_n\in\mathcal T}
       \|U_{h,M}^n-U_h(\vartheta_n)\|_{2,h},
 \qquad\mathcal T=\{\Theta/4,\Theta/2,\Theta\},
 \label{eq:p3checks-time-memory-errors}
\end{align}
where $\Theta=H(1)=1.0515132504186746$. These are absolute errors
maximized over three prescribed observations, each an actual time node.
The continuous-time references use numerical inverse Laplace transforms
of the finite-compartment resolvent and the direct resolvent difference
for continuous memory. Thus a refined time-stepping solution is not used
as the reference for this temporal test.

For $Q_h=P g_h$ with the prescribed delayed pulse, exact source cell
averages satisfy $\|g_h\|_{2,h}\le C_g=3.15809406831096\ldots$.
The source factors in Section~\ref{sec:p3ce-combined} give
\begin{align}
 B_{\rm FE}(\tau,M)
 &=C_2\varepsilon_M^{(1)}(\Theta)+\tfrac12 C_2\tau,
 &C_2&=C_g\,200\pi(\Theta-0.15),\nonumber\\
 B_{\rm BDF2}(\tau,M)
 &=C_2\varepsilon_M^{(1)}(\Theta)+C_B\tau^2,
 &C_B&=\tfrac{\sqrt5}{4}C_g\,3000\pi^2(\Theta-0.15).
 \label{eq:p3checks-time-memory-bounds}
\end{align}
Numerically, $C_2=1788.8628485809\ldots$ and
$C_B=47124.1128394591\ldots$.
The Euler steps satisfy the explicit restriction. Every tested
$\tau$ is below $0.05$, so the first backward-Euler step precedes
source activation and its startup error vanishes for this pulse.
This simplification does not apply to the manufactured source, whose
startup contribution is retained separately.
The primitive-kernel term uses Eq.~\eqref{eq:p3ce-kernel-enclosure}.

For each $M=8,16,32,64$, Euler uses
$N_\vartheta=1112,2224,4448,8896,17792,35584$, while BDF2 uses
$N_\vartheta=100,200,400,800,1600,3200,6400$.
The study contains $52$ solver runs.
Table~\ref{tab:p3checks-time-memory} reports the finest time resolution
and the observed order
$p=\log_2(E_{\rm time}(2\tau)/E_{\rm time}(\tau))$.
The two methods approach their first- and second-order temporal rates.
The counts are selected for refinement, not a matched-cost comparison.

\begin{table}[htbp]
\centering\small
\caption{Pulse-source time--memory check on $100$ cells.
Errors are maxima over $\mathcal T=\{\Theta/4,\Theta/2,\Theta\}$ and
are reported to five significant digits. The last column evaluates
Eq.~\eqref{eq:p3checks-time-memory-bounds}; it is an upper estimate,
not a measured error.}
\label{tab:p3checks-time-memory}
\begin{tabular}{lrrcrrr}
\toprule
Method & $M$ & $N_\vartheta$ & $p$ & $E_{\rm time}$ & $E_{\rm tm}$ & $B_{\rm method}$\\
\midrule
Euler & 8 & 35584 & 1.0023 & $9.6071\times10^{-4}$ & $6.6217\times10^{-4}$ & $1.0535\times10^{-1}$ \\
Euler & 16 & 35584 & 1.0023 & $9.6389\times10^{-4}$ & $9.9425\times10^{-4}$ & $2.9766\times10^{-2}$ \\
Euler & 32 & 35584 & 1.0023 & $9.6377\times10^{-4}$ & $9.6376\times10^{-4}$ & $2.6546\times10^{-2}$ \\
Euler & 64 & 35584 & 1.0023 & $9.6377\times10^{-4}$ & $9.6377\times10^{-4}$ & $2.6434\times10^{-2}$ \\
\midrule
BDF2 & 8 & 6400 & 1.9992 & $1.6266\times10^{-5}$ & $1.2013\times10^{-3}$ & $8.0187\times10^{-2}$ \\
BDF2 & 16 & 6400 & 1.9992 & $1.6275\times10^{-5}$ & $5.5156\times10^{-5}$ & $4.6073\times10^{-3}$ \\
BDF2 & 32 & 6400 & 1.9992 & $1.6275\times10^{-5}$ & $1.6264\times10^{-5}$ & $1.3871\times10^{-3}$ \\
BDF2 & 64 & 6400 & 1.9992 & $1.6275\times10^{-5}$ & $1.6275\times10^{-5}$ & $1.2758\times10^{-3}$ \\
\bottomrule
\end{tabular}
\end{table}

For $M=8$, increasing the BDF2 count from $800$ to $6400$ reduces
$E_{\rm time}$ from $1.0359\times10^{-3}$ to
$1.6266\times10^{-5}$, whereas $E_{\rm tm}$ approaches
$1.2013\times10^{-3}$. The independently computed sampled pure-memory
discrepancy is $1.2065\times10^{-3}$, explaining this plateau.
At $M=16$ the finest combined discrepancy is $5.5156\times10^{-5}$;
for $M=32$ and $64$ it is approximately $1.63\times10^{-5}$ and
tracks the remaining temporal error. At intermediate resolutions,
component errors can partly cancel, so combined error need not decrease
monotonically. In particular, the smaller combined Euler error at
$M=8$ does not establish superiority of the coarser memory approximation.

\begin{figure}[htbp]
\centering
\includegraphics[width=\linewidth]{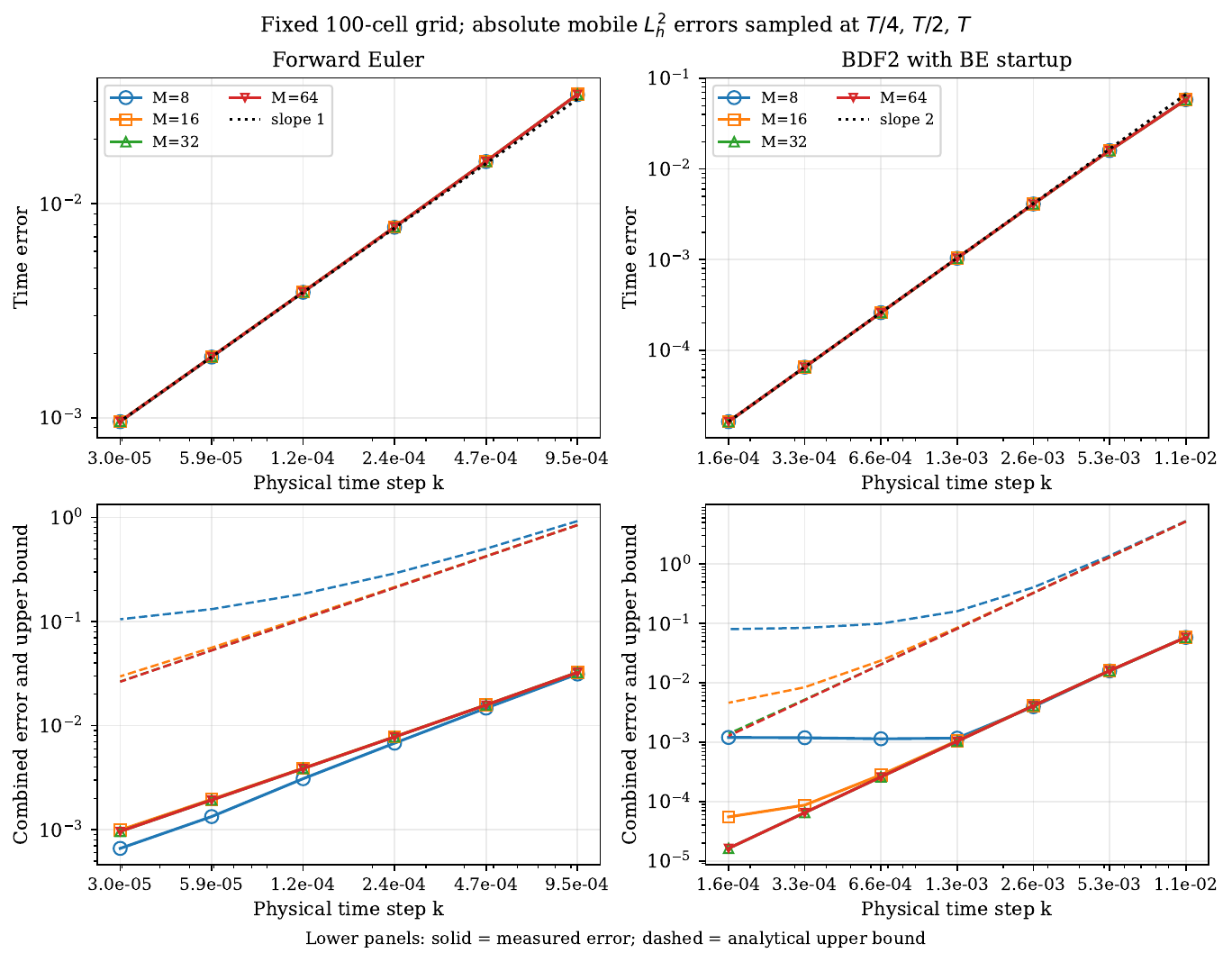}
\caption{Time and memory verification for the pulse source on
a fixed $100$-cell grid, with $M=8,16,32,64$.
Top: temporal errors against the corresponding continuous-time
finite-compartment state; dotted guides indicate slopes one and two.
Bottom: time--memory errors against the continuous-rate state on the
same mesh (solid) and the evaluated analytical upper estimates
(dashed). Each measured value is a maximum over three observation
times. The finest counts are $35584$ for Euler and $6400$ for BDF2.
Several curves overlap because their temporal errors are nearly equal.
The figure notation $k$ and $T$ denotes $\tau$ and $\Theta$, respectively.}
\label{fig:p3checks-time-memory}
\end{figure}

Figure~\ref{fig:p3checks-time-memory} displays all refinements, including $N_\vartheta=6400$ for BDF2. Every reported
sampled maximum lies below its corresponding evaluated time or
time--memory estimate. The estimates remain conservative. Numerical
agreement at these three observations is distinct from an
interval-certified check at all time nodes.

\subsection{Manufactured verification with continuous memory}
\label{sec:p3full-manufactured}
We verify the transport discretization with memory retained, using the
application coefficients and the physical-time interval
$0\leq\vartheta\leq\vartheta_f=H(1)=1.0515132504186746$.
Write $Lw=-\partial_x(vw-Dw_x)$ and prescribe
\begin{equation}
 u_\star(x,\vartheta)=\vartheta^2\varphi(x),\qquad
 \varphi(x)=\sin(\pi x)+x(1-x).
 \label{eq:p3full-exact}
\end{equation}
The mobile and immobile initial states are zero, and the mobile boundary
values are homogeneous. In particular, the trapping parameters remain
$\kappa_m=0.8$, $\lambda_0=0.5$ and $\gamma=0.5$.
For the continuous-memory equation
$u_\vartheta=Lu-\kappa_du-\ell*u_\vartheta+q$, the manufactured source is
\begin{equation}
 q_\star(x,\vartheta)
 =\bigl[2\vartheta+\kappa_d\vartheta^2+q_\ell(\vartheta)\bigr]
       \varphi(x)-\vartheta^2L\varphi(x),\qquad
 q_\ell(\vartheta)=2\int_0^\vartheta
       \ell(r)(\vartheta-r)\,dr.
 \label{eq:p3full-forcing}
\end{equation}
The scalar integral is evaluated through its incomplete-gamma expression.
The source uses the continuous kernel for every compartment count.

Let $\Pi_h$ denote cell averaging, $\Phi_h=\Pi_h\varphi$, and
$F_\varphi=v\varphi-D\varphi'$.
For a cell with faces $x_{i-1/2}$ and $x_{i+1/2}$, its source average is
\begin{equation}
 Q_{\star,i}(\vartheta)
 =\bigl[2\vartheta+\kappa_d\vartheta^2+q_\ell(\vartheta)\bigr](\Phi_h)_i
 +\frac{\vartheta^2}{h}
   \bigl[F_\varphi(x_{i+1/2})-F_\varphi(x_{i-1/2})\bigr].
 \label{eq:p3full-source-average}
\end{equation}
Thus the forcing is independent of the discrete transport matrix and
does not remove spatial or compartment errors by construction.
It can be signed, so the nonnegative-source positivity result is not
invoked for this benchmark. Although $Q_\star(0)=0$, one has
$Q_\star'(0)=2\Phi_h$; the backward-Euler initialization therefore remains
part of the temporal-error assessment.

\paragraph*{Reference fields and error definitions.}\enspace 
Let $U_h(\vartheta)$ denote the continuous-time spatial discretization
with continuous memory, and let $U_{h,M}(\vartheta)$ denote the
continuous-time $M$-compartment approximation on the same mesh.
The computed time-stepping state is $U_{h,M}^n$ at
$\vartheta_n=n\tau$. With
$\|V\|_{2,h}=(h\sum_i|V_i|^2)^{1/2}$, define
\begin{align}
 E_{\rm sp}(\vartheta)
   &=\|U_h(\vartheta)-\Pi_hu_\star(\vartheta)\|_{2,h},\nonumber\\
 E_{\rm mem}(\vartheta;M)
   &=\|U_{h,M}(\vartheta)-U_h(\vartheta)\|_{2,h},\nonumber\\
 E_{\rm time}^n
   &=\|U_{h,M}^n-U_{h,M}(\vartheta_n)\|_{2,h},\nonumber\\
 E_{\rm total}^n
   &=\|U_{h,M}^n-\Pi_hu_\star(\vartheta_n)\|_{2,h}.
 \label{eq:p3full-error-components}
\end{align}
All errors reported in this subsection are absolute cell-average norms.
They differ from the relative space--time discrepancies used for the
application below. The component error vectors add to the total error
vector; their norms need not add, and cancellation is possible.

The two continuous-time fields are evaluated by numerical inverse
Laplace transformation, without time stepping. This calculation uses
the same spatial matrix and source averages as the full solvers.
To avoid subtracting nearly equal reference states, we invert the error
transform directly. Define
$\mathscr K(s)=s\widehat\ell(s)$,
$\mathscr K_M(s)=\sum_{j=1}^M b_j\lambda_js/(s+\lambda_j)$ and
$d_h=L_h\Phi_h-\Pi_hL\varphi$. Then
\begin{equation}
 \widehat{U_{h,M}-\Pi_hu_\star}(s)
 =\frac{2}{s^3}
  \bigl[(s+\kappa_d+\mathscr K_M(s))I-L_h\bigr]^{-1}
  \bigl[d_h+(\mathscr K(s)-\mathscr K_M(s))\Phi_h\bigr].
 \label{eq:p3full-reference-transform}
\end{equation}
Replacing $\mathscr K_M$ by $\mathscr K$ gives the continuous-memory
spatial reference. The time-and-compartment study uses degree-$88$
Talbot inversion~\cite{Talbot1979}, checked against degree $64$ at every observation time
and compartment resolution. Selected components are also checked by
de Hoog inversion~\cite{DeHoogKnightStokes1982}. An independent full coupled calculation with
$100$ cells and $8$ compartments, integrated by adaptive Radau and BDF,
agrees with the transformed reference to maximum absolute differences
of $8.89\times10^{-16}$ and $1.70\times10^{-14}$, respectively.
These checks concern numerical reference consistency; high-precision
inversion does not eliminate errors in floating-point coefficients or
certify the continuous PDE solution.

\paragraph*{Spatial refinement.}\enspace
We first retain continuous time and the continuous exchange spectrum,
so this test measures $E_{\rm sp}$ alone. The spatial mesh sequence is
$N_x=100,200,400,800,1600,3200,6400$.
Observations are recorded at physical times $0.25$, $0.50$ and
$\vartheta_f$; Table~\ref{tab:p3full-space} gives the final-time values
and $p_h=\log_2(E_{\rm sp}(2h)/E_{\rm sp}(h))$ for the row with
$h=1/N_x$.
The errors decrease by approximately a factor of four under mesh
halving, consistent with second-order cell-average convergence.

\begin{table}[htbp]
\centering\small
\caption{Spatial verification with continuous time and continuous
memory. Errors are absolute cell-average $L^2$ errors at
$\vartheta_f=1.0515132504186746$. The mesh sequence reaches
$6400$ cells.}
\label{tab:p3full-space}
\begin{tabular}{rrc}
\toprule
$N_x$ & $E_{\rm sp}(\vartheta_f)$ & $p_h$\\
\midrule
100 & $9.91557435\times10^{-5}$ & --\\
200 & $2.47701211\times10^{-5}$ & 2.0011\\
400 & $6.19163067\times10^{-6}$ & 2.0002\\
800 & $1.54785579\times10^{-6}$ & 2.0000\\
1600 & $3.86960780\times10^{-7}$ & 2.0000\\
3200 & $9.67399008\times10^{-8}$ & 2.0000\\
6400 & $2.41848837\times10^{-8}$ & 2.0000\\
\bottomrule
\end{tabular}
\end{table}

\begin{figure}[htbp]
\centering
\includegraphics[width=\linewidth]{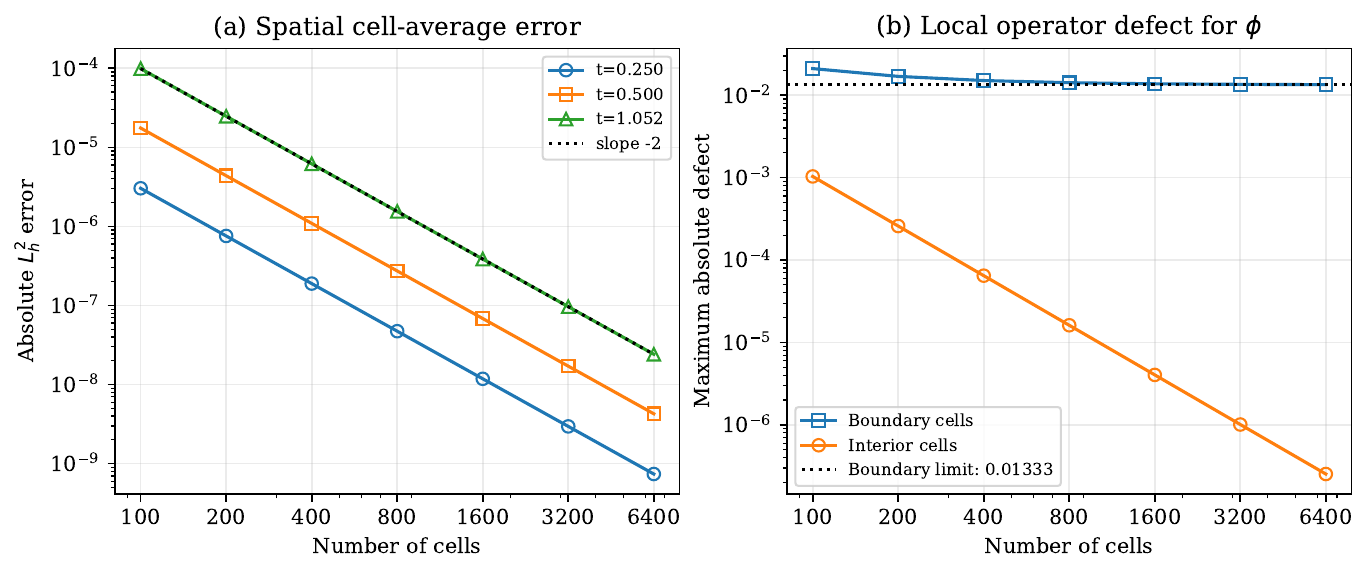}
\caption{Spatial verification with continuous time and continuous memory.
(a) Absolute cell-average errors at physical times
$\vartheta=0.25,0.50,\vartheta_f$ for $N_x=100$ through $6400$;
the dotted line has slope $-2$ against the cell count.
(b) Maximum absolute values of
$L_h\Pi_h\varphi-\Pi_hL\varphi$ on boundary and interior cells.
The boundary defect approaches $2D(0)/3$.
In the legend, $t$ denotes physical time $\vartheta$.}
\label{fig:p3checks-spatial-error}
\end{figure}

Figure~\ref{fig:p3checks-spatial-error} displays the full spatial
refinement sequence at the three observation times. The second panel
also checks the local consistency mechanism used in
Lemma~\ref{lem:p3ce-projection}: the interior cell defects decrease,
while the boundary cell defects approach a nonzero limit.
For this $\varphi$, the limit is $2D(0)/3=0.013\overline3$.
Thus a uniformly second-order local cell-residual estimate would be
incorrect, even though the global cell-average state error converges at
second order. The boundary lifting in the projection argument accounts
for this distinction.

\paragraph*{Temporal refinement.}\enspace 
Next fix $N_x=400$ and $M=64$. Forward Euler uses
$N_\vartheta=13148,26296,52592,105184$, while backward-Euler/BDF2 uses
$N_\vartheta=800,1600,3200,6400$.
The initial explicit count satisfies the positivity step restriction
with safety factor $0.8$; subsequent steps are obtained by halving.
Both methods are observed at the common times
$\mathcal T=\{\vartheta_f/4,\vartheta_f/2,\vartheta_f\}$.
These times differ from the earlier spatial-study observation times.
Table~\ref{tab:p3full-time} reports final-time errors against the
matching continuous-time $64$-compartment reference and
$p_\tau=\log_2(E_{\rm time}(2\tau)/E_{\rm time}(\tau))$ for the row
with step $\tau$. The measured orders approach one and two, respectively.

\begin{table}[htbp]
\centering\small
\caption{Temporal verification for the manufactured problem at
$N_x=400$ and $M=64$. Errors at $\vartheta_f$ use the matching
continuous-time compartment reference; spatial and memory discrepancies
are excluded from this column.}
\label{tab:p3full-time}
\begin{tabular}{lrrc}
\toprule
Method & $N_\vartheta$ & $E_{\rm time}(\vartheta_f)$ & $p_\tau$\\
\midrule
Explicit Euler & 13148 & $2.11849553\times10^{-5}$ & --\\
Explicit Euler & 26296 & $1.05924644\times10^{-5}$ & 1.0000\\
Explicit Euler & 52592 & $5.29622893\times10^{-6}$ & 1.0000\\
Explicit Euler & 105184 & $2.64811364\times10^{-6}$ & 1.0000\\
\midrule
Implicit BDF2 & 800 & $1.00174245\times10^{-7}$ & --\\
Implicit BDF2 & 1600 & $2.50303808\times10^{-8}$ & 2.0008\\
Implicit BDF2 & 3200 & $6.25602267\times10^{-9}$ & 2.0004\\
Implicit BDF2 & 6400 & $1.56379326\times10^{-9}$ & 2.0002\\
\bottomrule
\end{tabular}
\end{table}

\paragraph*{Compartment refinement and complete-solver errors.}\enspace 
At the same spatial resolution, we compare the continuous-time
$M$-compartment fields for $M=8,16,32,64$ with the continuous-memory
field. Table~\ref{tab:p3full-memory} reports final-time discrepancies
and the largest discrepancy on $\mathcal T$. The latter is a maximum
over three observations, not over the entire time interval.
For $M=64$, the discrepancy approaches the precision level of the
floating-point coefficients and is not evidence of corresponding
accuracy for the PDE solution. The rapid decrease in this smooth test
does not prescribe a sufficient compartment count for the pulse
application or for other initial-time behavior.

\begin{table}[htbp]
\centering\small
\caption{Continuous-time compartment refinement on $400$ cells.
Both columns compare with continuous memory on the same spatial mesh.
The observation set is
$\mathcal T=\{\vartheta_f/4,\vartheta_f/2,\vartheta_f\}$.}
\label{tab:p3full-memory}
\begin{tabular}{rrr}
\toprule
$M$ & $E_{\rm mem}(\vartheta_f;M)$
& $\max_{\vartheta\in\mathcal T}E_{\rm mem}(\vartheta;M)$\\
\midrule
8 & $8.64253658\times10^{-7}$ & $1.36011858\times10^{-6}$\\
16 & $1.53081459\times10^{-9}$ & $1.48197474\times10^{-8}$\\
32 & $4.60281823\times10^{-13}$ & $1.40913405\times10^{-11}$\\
64 & $3.29767440\times10^{-15}$ & $3.29767440\times10^{-15}$\\
\bottomrule
\end{tabular}
\end{table}

A full-solver compartment sweep additionally uses the finest time count
of each method, reusing the $M=64$ runs from the temporal study.
The two sweeps therefore require $14$ distinct fully discrete
trajectories. At $N_x=400$, $M=64$ and $\vartheta=\vartheta_f$, the
finest implicit run has total error $6.19214128\times10^{-6}$, compared
with spatial discrepancy $6.19163067\times10^{-6}$ and temporal error
$1.56379326\times10^{-9}$. Its total error is dominated by the spatial
component. The finest explicit total error is
$5.63683215\times10^{-6}$, although its temporal error is
$2.64811364\times10^{-6}$. The smaller total norm reflects partial
cancellation of component errors; it does not establish a general
accuracy advantage for the explicit method.
Across these $14$ runs the largest discrete mass-balance defect is
$3.14\times10^{-13}$. This accounting check is separate from state
accuracy.

\begin{figure}[htbp]
\centering
\includegraphics[width=\linewidth]{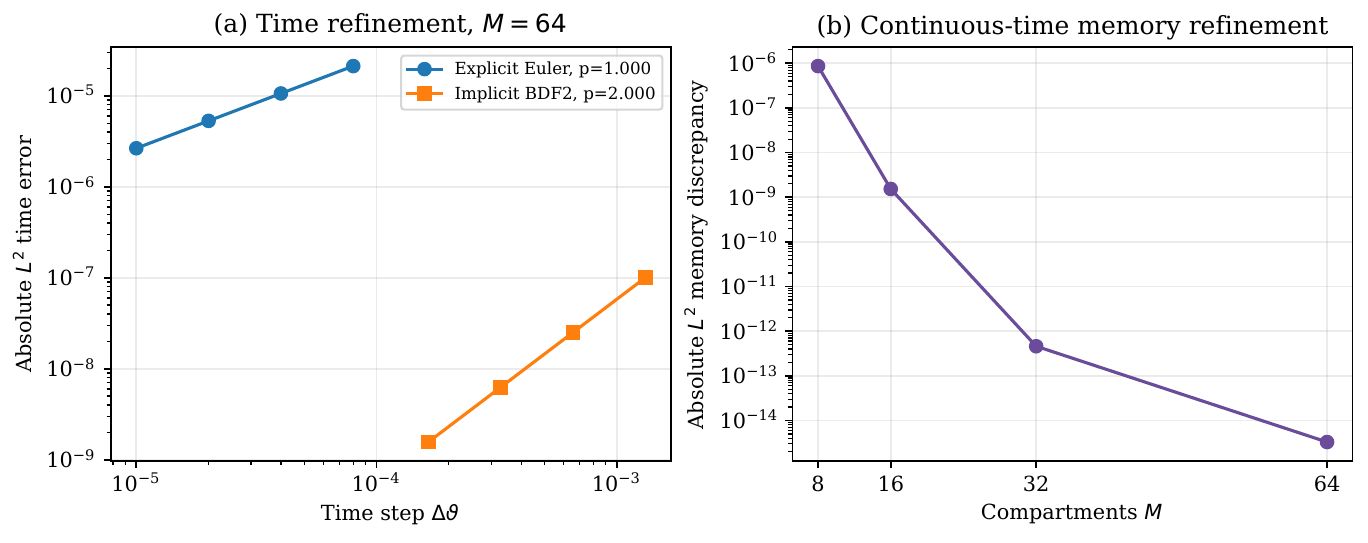}
\caption{Manufactured verification on $400$ cells at
$\vartheta_f=1.0515132504186746$.
(a) Temporal errors relative to the matching continuous-time
$64$-compartment reference; the displayed order uses the final two
time resolutions. (b) Continuous-time compartment discrepancies
relative to continuous memory on the same spatial mesh.
The smallest discrepancy approaches floating-point coefficient
precision. Neither panel displays total PDE error.}
\label{fig:p3full-refinement}
\end{figure}

These component studies show the expected refinement behavior for the
prescribed smooth state. We next compare the full error with the
combined analytical estimate.

\paragraph*{Evaluation of the combined estimate.}\enspace 
We evaluate Theorem~\ref{thm:p3ce-combined} on the same computed trajectories,
using the manufactured constants in
Eqs.~\eqref{eq:p3ce-manufactured-c2}--\eqref{eq:p3ce-manufactured-time}.
For the prescribed diffusion, the analytic bounds
$d_*=0.02\exp(-256/243)$ and $d^*=0.02$ are valid.
At $N_x=400$, these give the spatial contribution
$h^2C_{\rm sp}=2.55288046\times10^{-2}$.
The source factors are $\overline C_{2,h}=8.73187137$ and
$\overline C_{3,h}=5.60729928$.
The primitive-kernel defect is evaluated with piecewise-linear
interpolation and the remainder control described in
Eq.~\eqref{eq:p3ce-kernel-enclosure}.
The resulting memory contributions decrease from
$3.8520\times10^{-4}$ at $M=8$ to $1.8399\times10^{-8}$ at $M=64$.

For the manufactured state \eqref{eq:p3full-exact}, set
$\Theta=\vartheta_f$ and define
\begin{equation}
	\varepsilon_{\varphi,h}
	:=\|R_h\varphi-\Pi_h\varphi\|_{2,h}.
	\label{eq:p3rev-projection-error}
\end{equation}
Linearity gives
$\eta(\vartheta):=R_hu_\star-\Pi_hu_\star
=\vartheta^2(R_h\varphi-\Pi_h\varphi)$.
The projected residual is
$\eta'+\kappa_d\eta+\ell*\eta'$.
Applying \eqref{eq:p3ce-forcing} to the zero-initial-value
error $U_h-R_hu_\star$, followed by the triangle inequality
and $\|\ell\|_{L^1(0,\Theta)}\leq B$, yields
\begin{equation}
	\begin{aligned}
		\max_{0\leq\vartheta\leq\Theta}E_{\rm sp}(\vartheta)
		&\leq
		\left[\Theta^2+(1+B)\int_0^\Theta 2s\,ds
		+\kappa_d\int_0^\Theta s^2\,ds\right]
		\varepsilon_{\varphi,h}\\
		&=\left[(2+B)\Theta^2+\frac{\kappa_d\Theta^3}{3}\right]
		\varepsilon_{\varphi,h}.
	\end{aligned}
	\label{eq:p3rev-projection-spatial}
\end{equation}
The bound $\varepsilon_{\varphi,h}\leq h^2A(\varphi)$
recovers $h^2C_{\rm sp}(u_\star,\Theta)$.
Direct evaluation of $\varepsilon_{\varphi,h}$ provides
a benchmark-specific spatial contribution; the memory
and temporal terms remain unchanged.

I run the code. Here is Table 8. \begin{table}[htbp]
	\centering\small
	\caption{Combined error estimates on $400$ cells.
		Observed errors are evaluated at $\vartheta=\Theta$.
		The analytical estimate uses $h^2C_{\rm sp}$; the projection
		estimate uses \eqref{eq:p3rev-projection-spatial}.
		Both retain the same memory and temporal terms and estimate
		the maximum error over the time nodes.}
	\label{tab:p3full-combined-bound}
	\begin{tabular}{lrrrrr}
		\toprule
		Method & $M$ & $N_\vartheta$ & Observed error
		& \shortstack{Analytical\\estimate}
		& \shortstack{Projection\\estimate}\\
		\midrule
		Euler & 8 & 105184 & $5.4731\times10^{-6}$ & $2.5958\times10^{-2}$ & $4.6046\times10^{-4}$ \\
		Euler & 16 & 105184 & $5.6355\times10^{-6}$ & $2.5589\times10^{-2}$ & $9.1543\times10^{-5}$ \\
		Euler & 32 & 105184 & $5.6368\times10^{-6}$ & $2.5573\times10^{-2}$ & $7.5824\times10^{-5}$ \\
		Euler & 64 & 105184 & $5.6368\times10^{-6}$ & $2.5572\times10^{-2}$ & $7.5281\times10^{-5}$ \\
		\midrule
		BDF2 & 8 & 6400 & $6.3972\times10^{-6}$ & $2.5914\times10^{-2}$ & $4.1696\times10^{-4}$ \\
		BDF2 & 16 & 6400 & $6.1910\times10^{-6}$ & $2.5545\times10^{-2}$ & $4.8036\times10^{-5}$ \\
		BDF2 & 32 & 6400 & $6.1921\times10^{-6}$ & $2.5530\times10^{-2}$ & $3.2317\times10^{-5}$ \\
		BDF2 & 64 & 6400 & $6.1921\times10^{-6}$ & $2.5529\times10^{-2}$ & $3.1774\times10^{-5}$ \\
		\bottomrule
	\end{tabular}
\end{table}
The analytical estimate is dominated by its spatial contribution.
Replacing $h^2A(\varphi)$ with the computed projection error
$\varepsilon_{\varphi,h}$ substantially reduces this conservatism.
At $M=64$, the combined estimates decrease by factors of approximately
$340$ for Euler and $803$ for BDF2. The revised estimates are
$13.4$ and $5.13$ times the respective final-time errors.
The memory and temporal contributions remain unchanged.
These results demonstrate a tighter estimate for the manufactured
benchmark, without establishing sharpness.
\subsection{Separate application refinements}
\label{sec:p3pair-refinement}
For the pulse application, the implicit time study uses $6400$, $12800$
and $25600$ steps at $1600$ cells and $64$ compartments.
At this spatial and compartment resolution,
$\tau_{\rm pos}=6.4432484\times10^{-6}$.
A safety factor of $0.9$ gives the initial explicit count $181330$;
the subsequent counts are $362660$ and $725320$.
Space and compartment counts are then varied separately, holding the
other resolutions at their method-specific finest values.

For field comparisons, both numerical solutions are interpolated
bilinearly onto a common uniform $257\times257$ grid in
$[0,1]\times[0,\vartheta_f]$. Relative $L^2$ discrepancies use tensor
trapezoidal integration in $dx\,d\vartheta$; $E_{\max}$ is the maximum
absolute difference on that grid. Each row in
Table~\ref{tab:p3pair-refinement} is compared with the finest run of its
own method.

\begin{table}[htbp]
\centering\small
\caption{Application refinement against method-specific references:
$(N_x,N_\vartheta,N_c)=(1600,25600,64)$ for the implicit method and
$(1600,725320,64)$ for the explicit method.
Entries are numerical discrepancies, not certified errors.}
\label{tab:p3pair-refinement}
\begin{tabular}{llrrrcc}
\toprule
Method & Varied & $N_x$ & $N_\vartheta$ & $N_c$ & Relative $L^2$ & $E_{\max}$\\
\midrule
Implicit & Time & 1600 & 6400 & 64 & $1.1417\times10^{-5}$ & $4.2370\times10^{-5}$ \\
Implicit & Time & 1600 & 12800 & 64 & $2.2833\times10^{-6}$ & $8.4737\times10^{-6}$ \\
Implicit & Space & 400 & 25600 & 64 & $4.9162\times10^{-4}$ & $1.0220\times10^{-2}$ \\
Implicit & Space & 800 & 25600 & 64 & $9.1201\times10^{-5}$ & $1.4207\times10^{-3}$ \\
Implicit & Rates & 1600 & 25600 & 16 & $3.3005\times10^{-5}$ & $1.1188\times10^{-4}$ \\
Implicit & Rates & 1600 & 25600 & 32 & $1.3099\times10^{-7}$ & $5.1116\times10^{-7}$ \\
\midrule
Explicit & Time & 1600 & 181330 & 64 & $1.0717\times10^{-4}$ & $4.2679\times10^{-4}$ \\
Explicit & Time & 1600 & 362660 & 64 & $3.5720\times10^{-5}$ & $1.4225\times10^{-4}$ \\
Explicit & Space & 400 & 725320 & 64 & $4.9165\times10^{-4}$ & $1.0220\times10^{-2}$ \\
Explicit & Space & 800 & 725320 & 64 & $9.1208\times10^{-5}$ & $1.4208\times10^{-3}$ \\
Explicit & Rates & 1600 & 725320 & 16 & $3.3012\times10^{-5}$ & $1.1191\times10^{-4}$ \\
Explicit & Rates & 1600 & 725320 & 32 & $1.3103\times10^{-7}$ & $5.1135\times10^{-7}$ \\
\bottomrule
\end{tabular}
\end{table}

Successive absolute $L^2$ time-refinement gaps give observed orders
$2.0000$ and $1.0003$ for the implicit and explicit methods,
respectively. The closest explicit temporal gap is about $39\%$ of
its closest spatial gap; the corresponding implicit ratio is about
$2.5\%$. For both methods the temporal gap exceeds the very small
$32$-to-$64$ compartment discrepancy. Thus the smallest compartment
change does not demonstrate that the total error is equally small.

The finest explicit state differs from the finest implicit state by
$3.5725\times10^{-5}$ in relative $L^2$, with maximum absolute difference
$1.4226\times10^{-4}$. The maximum detector difference is
$1.3571\times10^{-4}$; final immobile and captured masses differ by
$9.5227\times10^{-8}$ and $8.8571\times10^{-8}$.
Across the application refinement and control runs, the largest
discrete balance defects are $3.55\times10^{-12}$ for the implicit
method and $2.91\times10^{-14}$ for the explicit method.
The monitored minimum concentration is zero in both families.

For the finest paired calculations, the recorded process CPU times are
$14.50$~s for the implicit method and $242.81$~s for the explicit method.
These measurements include numerical setup, time stepping, diagnostics
and output-array assembly, and exclude file writing and plotting.
The runs use different time-step counts and are not matched at a common
error tolerance; these timings therefore do not constitute an
equal-accuracy efficiency comparison.

\begin{figure}[htbp]
\centering
\includegraphics[width=\linewidth]{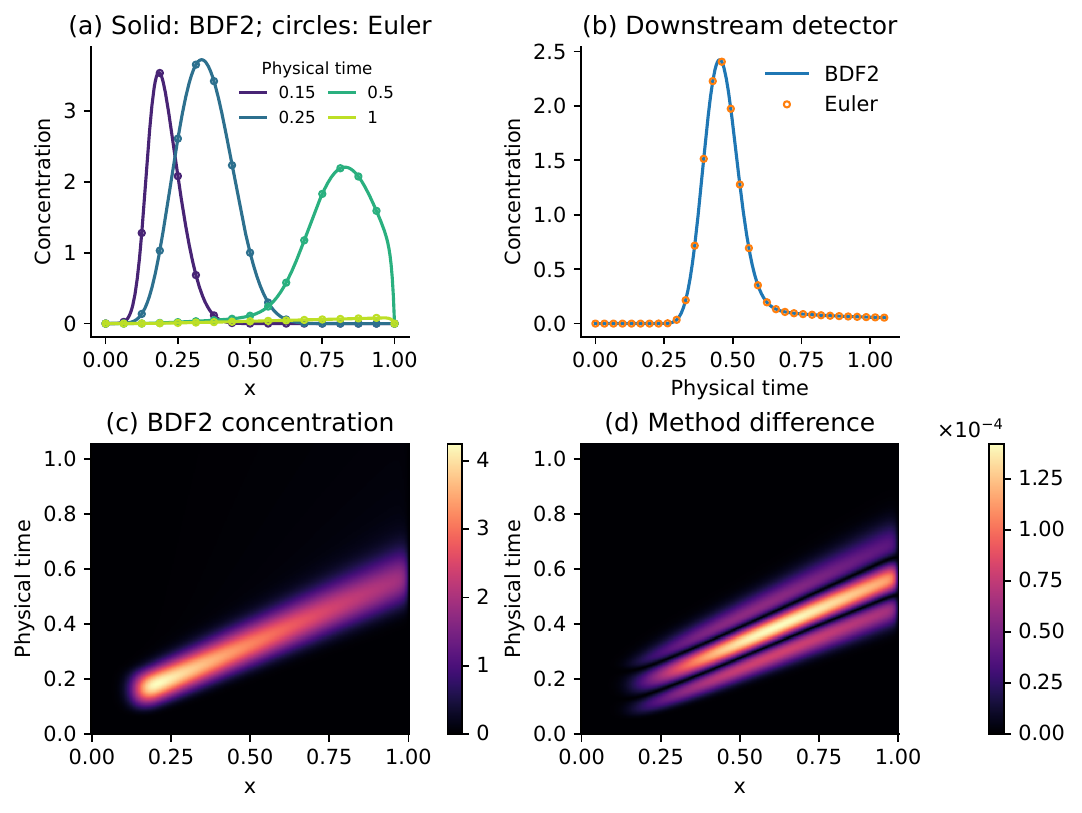}
\caption{Comparison of the finest implicit and explicit application
calculations, with $1600$ cells and $64$ compartments.
The implicit calculation uses $25600$ steps and the explicit calculation
uses $725320$. (a) Concentration profiles at physical times
$0.15$, $0.25$, $0.50$ and $1.00$: solid curves are implicit and circles
are explicit. (b) Downstream detector concentration.
(c) Implicit concentration in space--time.
(d) Absolute difference between the two numerical fields.
Panel (d) is a method discrepancy, not an exact-error map.}
\label{fig:p3pair-comparison}
\end{figure}

\subsection{Resolution checks for the application-output study}
\label{sec:app-reference-verification}
The physical-output study uses an implicit reference with
$(N_x,N_\vartheta,N_c)=(1600,6400,64)$ and
$\vartheta_f=H(1)=1.05151325$.
This reference is distinct from the $25600$-step implicit calculation
in Section~\ref{sec:p3pair-refinement}.
The checks and figures in this and the next subsection use the
$6400$-step reference. Space, time and compartments are refined separately.
For a coarser solution $U$ and the finest reference $U_f$, define
\begin{equation}
	E_{2,\mathrm{ref}}=
	\frac{\|U-U_f\|_{L^2(dx\,d\vartheta)}}
	{\|U_f\|_{L^2(dx\,d\vartheta)}},\qquad
	E_{\infty,\mathrm{ref}}=
	\max_{(x,\vartheta)\in\mathcal G}|U-U_f|,
	\label{eq:app-reference-metrics}
\end{equation}
where $\mathcal G$ is a uniform $257\times257$ grid on
$[0,1]\times[0,\vartheta_f]$. Both fields are interpolated bilinearly
onto this grid, and the integrals use the tensor trapezoidal rule.
This state norm uses physical time
$d\vartheta=dt/\psi(t)$.

\begin{table}[!htbp]
	\centering 
	\small
	\caption{Refinement of the implicit reference for the physical-output study. Each row is compared with
	the same finest calculation $(N_x,N_\vartheta,N_c)=(1600,6400,64)$.
The entries are refinement discrepancies, not certified errors.}
	\label{tab:app-refinement}
	\begin{tabular}{lrrrrr}
		\toprule
	Varied resolution & $N_x$ & $N_\vartheta$ & $N_c$
		& $E_{2,\mathrm{ref}}$ & $E_{\infty,\mathrm{ref}}$ \\
	\midrule
	Space & 400 & 6400 & 64 & $4.9162\times10^{-4}$ & $1.0220\times10^{-2}$ \\
	Space & 800 & 6400 & 64 & $9.1201\times10^{-5}$ & $1.4207\times10^{-3}$ \\
	Time & 1600 & 1600 & 64 & $1.8546\times10^{-4}$ & $7.0075\times10^{-4}$ \\
	Time & 1600 & 3200 & 64 & $3.7131\times10^{-5}$ & $1.4134\times10^{-4}$ \\
	Compartments & 1600 & 6400 & 16 & $3.3005\times10^{-5}$ & $1.1189\times10^{-4}$ \\
	Compartments & 1600 & 6400 & 32 & $1.3099\times10^{-7}$ & $5.1118\times10^{-7}$ \\
	\bottomrule
	\end{tabular}
\end{table}

Table~\ref{tab:app-refinement} shows decreasing discrepancies under
each refinement. The closest spatial and temporal comparisons give
relative differences of $9.12\times10^{-5}$ and $3.71\times10^{-5}$,
whereas increasing the compartment count from $32$ to $64$ changes
the state by $1.31\times10^{-7}$ in the same metric.
The capacity quadrature reproduces the total capacity, but not the
initial exchange coefficient exactly. For $64$ compartments,
$\ell_{64}(\vartheta)=\sum_j b_j\lambda_j e^{-\lambda_j\vartheta}$
differs from $\ell$ by $8.77\times10^{-3}$ at $\vartheta=0$ and
$3.48\times10^{-5}$ at $\vartheta=0.01$. Thus total-capacity agreement
alone is insufficient to assess the exchange approximation.

Reference reliability is also checked for the physical outputs.
Let $\delta_{\rm det}$ be the maximum absolute detector discrepancy on
the $257$ physical-time grid points, and let $\delta_{M_i}$ and
$\delta_{\rm cap}$ denote absolute discrepancies in final immobile and
captured mass. Table~\ref{tab:app-output-refinement} reports the closest
refinement in each direction, with all other resolutions fixed at their
finest values. These checks address observables that are not bounded
by the continuous $L^2$ state estimate alone.

\begin{table}[htbp]
	\centering
	\caption{Refinement discrepancies in physical outputs, measured against
		the $(1600,6400,64)$ reference. All entries are absolute differences.}
	\label{tab:app-output-refinement}
	\begin{tabular}{lccc}
\toprule
Comparison & $\delta_{\rm det}$ & $\delta_{M_i}$ & $\delta_{\rm cap}$ \\
	\midrule
$800$ versus $1600$ cells
& $2.1253\times10^{-4}$ & $6.4005\times10^{-6}$ & $9.8661\times10^{-6}$ \\
$3200$ versus $6400$ steps
& $1.2625\times10^{-4}$ & $4.7034\times10^{-9}$ & $1.2644\times10^{-10}$ \\
$32$ versus $64$ compartments
& $3.5528\times10^{-7}$ & $8.5764\times10^{-11}$ & $1.2585\times10^{-9}$ \\
		\bottomrule
	\end{tabular}
\end{table}

\begin{figure}[htbp]
	\centering
	\includegraphics[width=\linewidth]{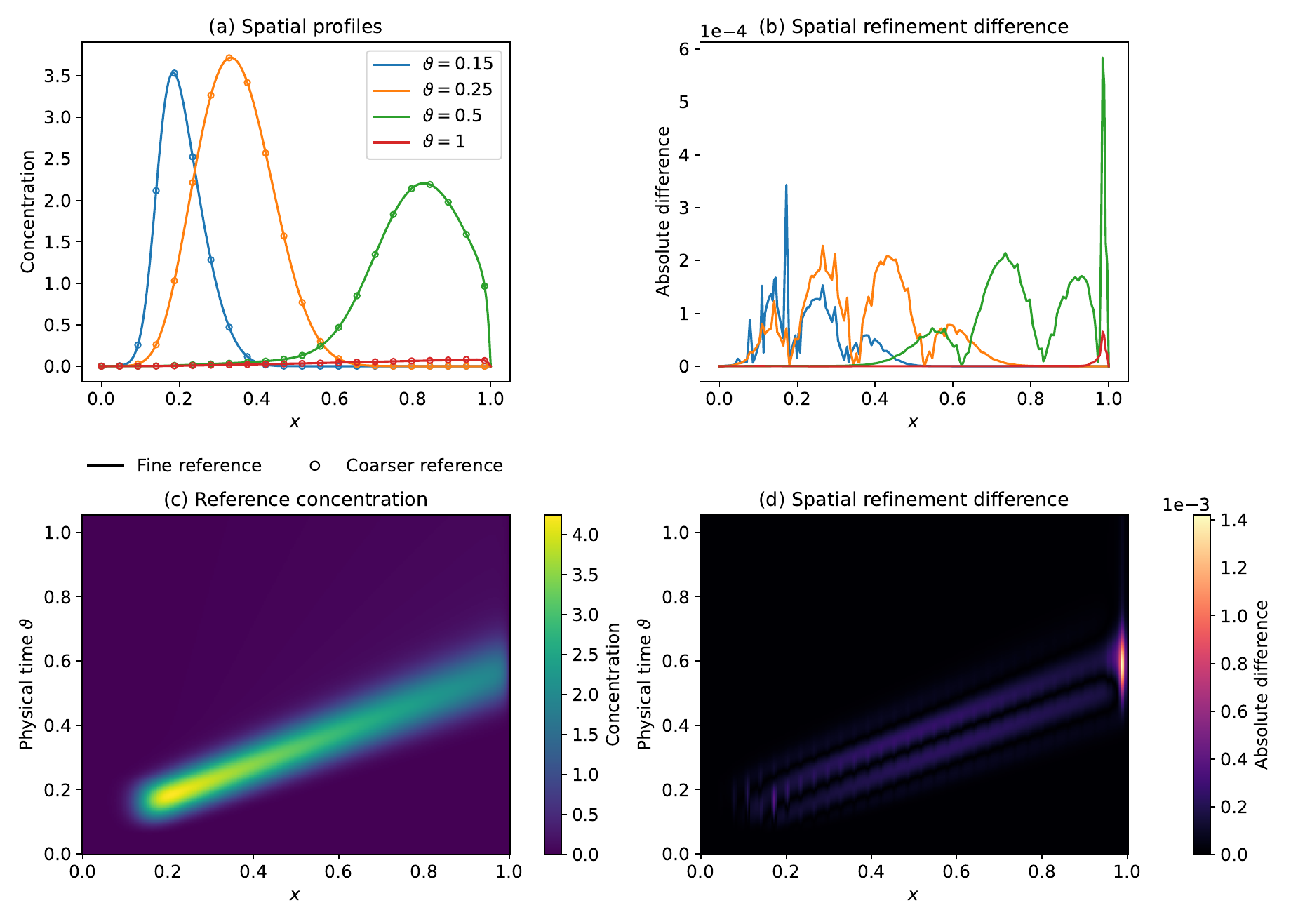}
	\caption{Deterministic reference for the porous-medium application.
	(a) Concentration profiles at physical times
	$\vartheta=0.15,0.25,0.50,1.00$: solid curves use $1600$ cells and
	circles use $800$ cells. (b) Absolute spatial-refinement differences
	at the same times. (c) Fine-reference concentration.
	(d) Absolute difference between the two spatial resolutions on a
	uniform $257\times257$ space--time grid. Both calculations use
	$6400$ time steps and $64$ compartments.}
	\label{fig:app-reference}
\end{figure}

Figure~\ref{fig:app-reference} shows the computed concentration and
the spatial-refinement differences between $800$ and $1600$ cells.
The discrepancy panels compare two deterministic approximations;
no exact application state is available.

\subsection{Transport response, storage and release}
\label{sec:app-transport-results}
Table~\ref{tab:app-observables} quantifies the trapping effect using
the $6400$-step implicit reference and its no-trapping control. The detector peak is
attenuated by $24.01\%$, while both sampled peak times remain close to
$0.452$. Their difference is one time step, so the calculations support
peak attenuation and persistent tailing rather than a resolved delay in
peak arrival. At the final time, $20.97\%$ of the injected mass remains
in the mobile and immobile phases. Cumulative capture is lower by
$0.20558$, or $20.56$ percentage points of the unit injected mass.

\begin{table}[htbp]
	\centering
	\caption{Physical observables with and without trapping. Final values
		are evaluated at $\vartheta_f=1.05151325$. Both cases use $1600$ cells
		and $6400$ time steps; the trapping case uses $64$ compartments.}
	\label{tab:app-observables}
	\begin{tabular}{lcc}
		\toprule
		Quantity & With trapping & Without trapping \\
		\midrule
		Detector peak & $2.42339$ & $3.18898$ \\
		Sampled peak time & $0.452315$ & $0.452479$ \\
		Final detector concentration & $5.47196\times10^{-2}$ & $8.2158\times10^{-11}$ \\
		Final mobile mass & $0.0322191$ & $3.0863\times10^{-8}$ \\
		Final immobile mass & $0.1775243$ & $0$ \\
		Final captured mass & $0.7529027$ & $0.9584822$ \\
		Final degraded mass & $0.0373539$ & $0.0415178$ \\
		\bottomrule
	\end{tabular}
\end{table}

\begin{figure}[htbp]
	\centering
	\includegraphics[width=\linewidth]{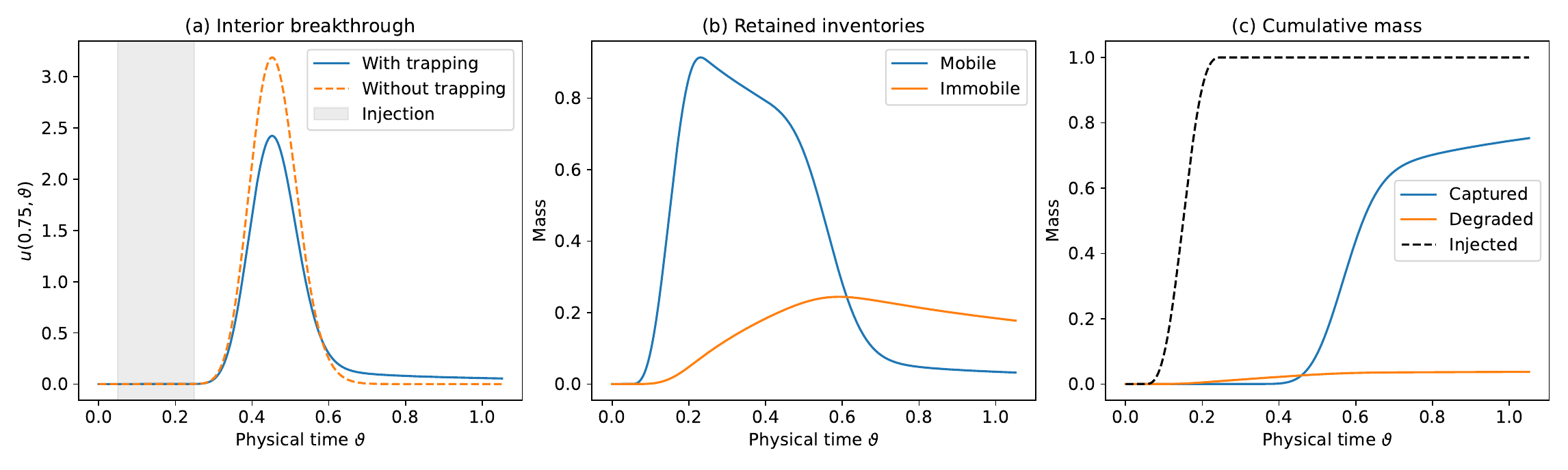}
	\caption{Reference transport observables. (a) Concentration at the
	interior detector $x=0.75$, with and without trapping; shading marks
	the injection interval. (b) Mobile and immobile mass for the trapping
	case. (c) Cumulative boundary capture, degradation and analytically
	integrated injection.}
	\label{fig:app-transport}
\end{figure}

Figure~\ref{fig:app-transport} connects the breakthrough curve to the
mass inventories. Concentration persists after the main pulse passes,
while immobile storage gradually decreases.
The minimum observed concentration is zero in both cases.

\begin{figure}[htbp]
	\centering
	\includegraphics[width=\linewidth]{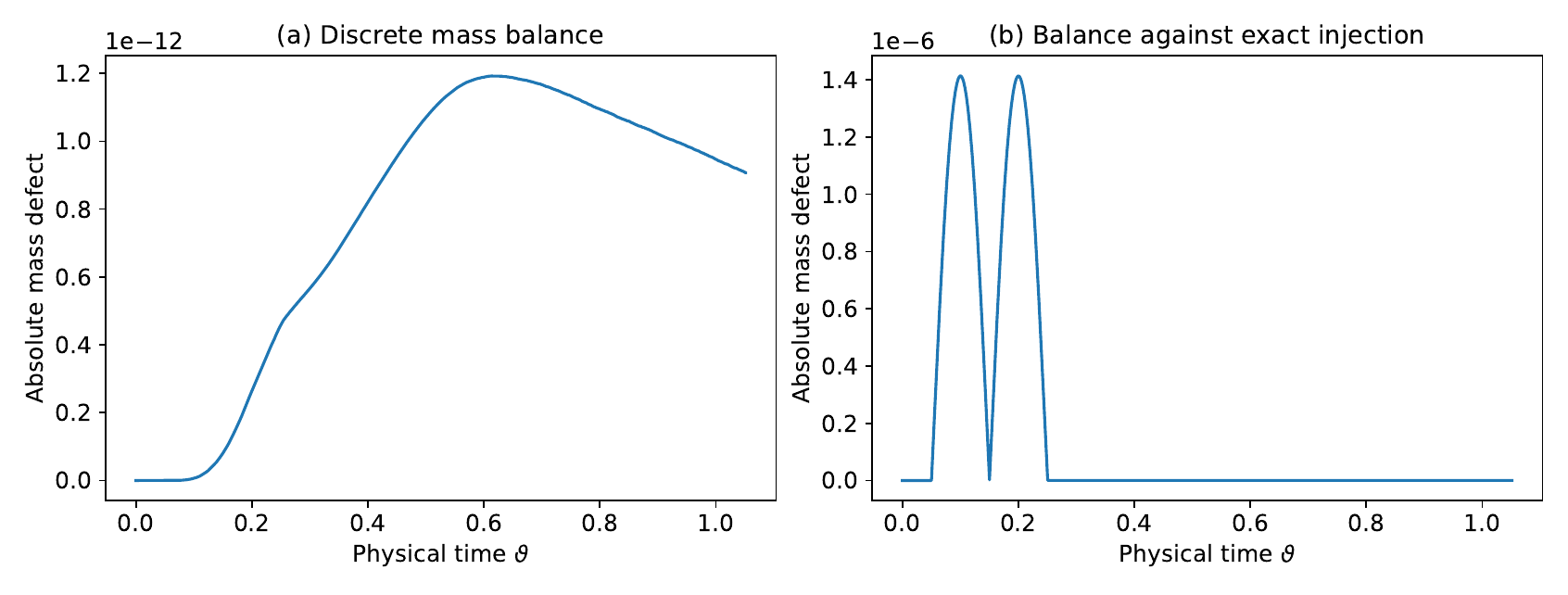}
	\caption{Conservation diagnostics for the $6400$-step implicit reference.
		(a) Absolute mass-balance defect with all cumulative terms integrated
		by the same discrete time formula. (b) Absolute balance defect using
		the exact cumulative injection. The respective maxima are
		$1.19\times10^{-12}$ and $1.41\times10^{-6}$.}
	\label{fig:app-balance}
\end{figure}

Two balance checks are distinguished in Fig.~\ref{fig:app-balance}.
When injection, capture and degradation use the same discrete time
formula as the states, the maximum absolute balance defect is
$1.19\times10^{-12}$. Replacing the discrete injection integral by
its analytical value gives a maximum defect of $1.41\times10^{-6}$.
The first check measures algebraic conservation of the discretization;
the second also exposes the time-integration discrepancy in the source.
Neither replaces the state-refinement study.

\subsubsection{Space--time structure of retention and release}
\begin{figure}[htbp]
	\centering
	\includegraphics[width=\linewidth]{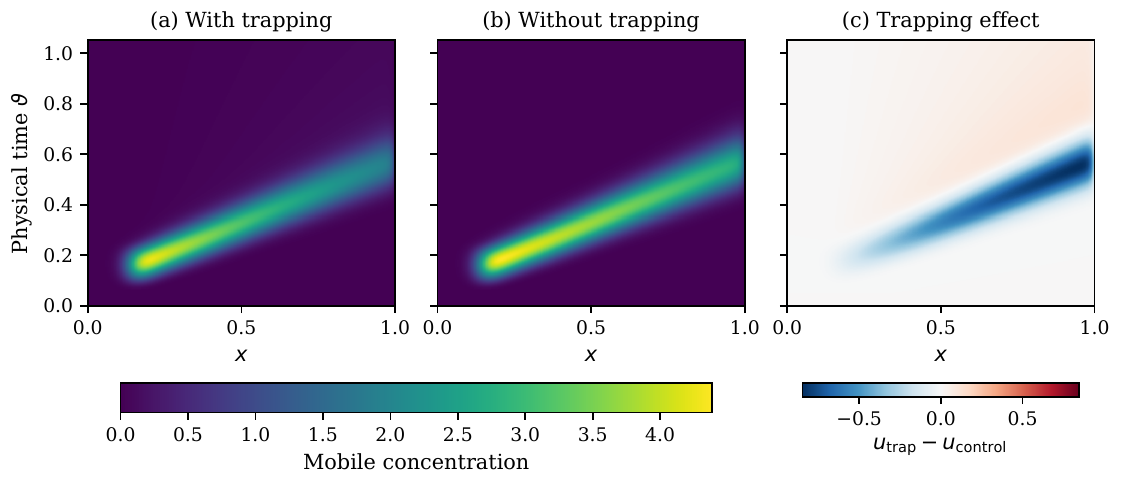}
	\caption{Effect of trapping in physical space--time.
	(a) Mobile concentration with trapping.
	(b) Mobile concentration without trapping, using the same color limits
	as panel (a). (c) The signed difference
	$\Delta u=u_{\rm trap}-u_{\rm control}$, with a symmetric linear color
	scale centered at zero. These deterministic-reference fields are
	displayed on the same $257\times257$ grid.}
	\label{fig:app-trapping-maps}
\end{figure}

The physical effect of trapping is resolved by comparing the mobile fields
$u_{\rm trap}$ and $u_{\rm control}$ with identical color limits in
Fig.~\ref{fig:app-trapping-maps}. Their signed difference,
$\Delta u=u_{\rm trap}-u_{\rm control}$, separates attenuation of the
moving concentration peak from the enhanced concentration in its wake.
This difference compares two prescribed physical models; it is not a
numerical-error estimate. All maps use position $x$ horizontally and
physical time $\vartheta$ vertically.

\begin{figure}[htbp]
	\centering
	\includegraphics[width=\linewidth]{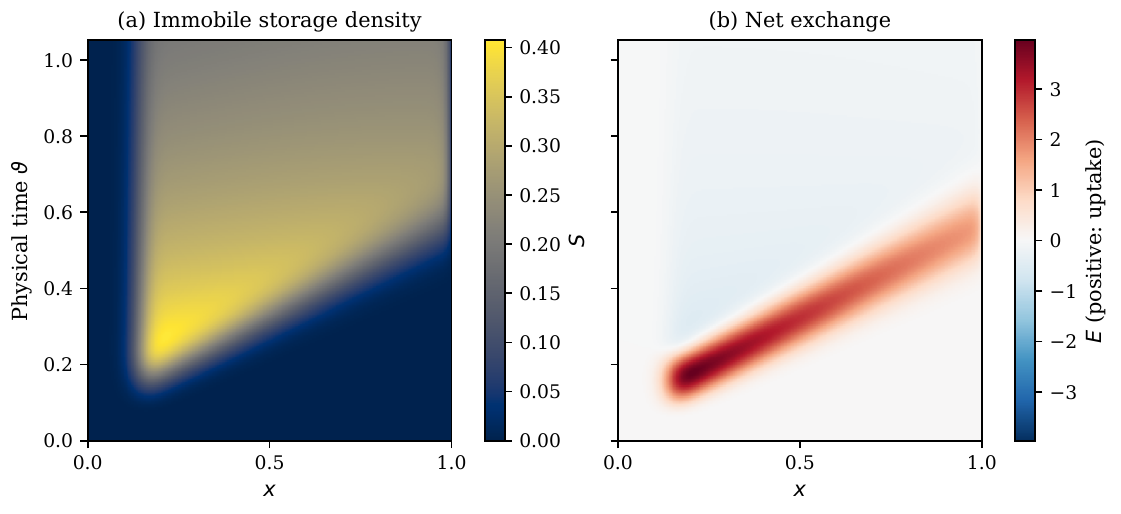}
	\caption{Immobile storage and exchange reconstructed from the computed
		reference trajectory. (a) Capacity-weighted immobile storage density
		$S$. (b) Net exchange $E$: positive values indicate uptake and negative
		values indicate release. The exchange color scale is linear and
		symmetric about zero. Both maps use physical time and $64$ compartments.}
	\label{fig:app-storage-exchange}
\end{figure}

To localize storage and release, define the compartment approximations
\begin{equation}
	S(x,\vartheta)=\sum_{j=1}^{N_c}b_jz_j(x,\vartheta),\qquad
	E(x,\vartheta)=\sum_{j=1}^{N_c}\lambda_jb_j
	[u(x,\vartheta)-z_j(x,\vartheta)].
	\label{eq:app-storage-exchange}
\end{equation}
Here $S$ is the capacity-weighted immobile storage density, whose spatial
integral gives the compartment approximation of $M_i$. Positive $E$ denotes uptake into immobile regions;
negative $E$ denotes release into the mobile phase.
The compartment states are reconstructed from the computed mobile trajectory
using the same rates, capacities and backward-Euler/BDF2 time formula as
the reference calculation. Their integrated storage agrees with the
saved immobile mass to within $1.97\times10^{-15}$.
Figure~\ref{fig:app-storage-exchange} shows that storage persists after
the mobile peak passes, with net release in the trailing region.

\FloatBarrier

\subsection{Scope of the numerical evidence}\label{sec:app-evidence}
The convergence studies and evaluated bounds apply to the specified
data, norms and observation sets. The manufactured state is known exactly; the application reference
is numerical.
Separate refinements assess detector traces, boundary capture and
immobile inventory. The mobile-state estimate alone does not certify
these observables.
Observed nonnegative BDF2 trajectories likewise do not establish
positivity for arbitrary data and steps.

The application is a nondimensional model study with prescribed
coefficients. Its results have not been validated against field data.
The heatmaps show one spatial coordinate and physical time.
They do not represent a two-dimensional spatial simulation.

\section{Conclusion}\label{sec:conclusion}
This work establishes an error analysis for conservative finite-volume
discretizations of reactive conformable transport with multirate memory.
The analysis separates spatial, compartment and temporal errors.
It accounts for the boundary flux closure and the backward-Euler
initialization of BDF2. For the two prescribed sources, the temporal
and memory prefactors remain uniform under mesh and compartment
refinement. The formulation also distinguishes the local conformable
time transformation from the hereditary response generated by
mobile--immobile exchange.

The refinement studies recover second-order spatial convergence,
first-order temporal convergence for Euler and second-order temporal
convergence for BDF2. They also demonstrate that temporal refinement
alone can leave a compartment error floor. Cancellation between error
components can temporarily conceal this limitation, which explains
the need to examine the components separately.

The manufactured benchmark identifies the spatial contribution as
the main source of conservatism in the analytical estimate.
On $400$ cells with $64$ compartments, evaluating the elliptic
projection error directly reduces the combined estimates by factors
of approximately $340$ for Euler and $803$ for BDF2.
The resulting estimates are approximately $13.4$ and $5.13$ times
the respective final-time errors. Thus, retaining the discrete
projection information substantially improves the quantitative
value of the estimate without changing the computed solution.

In the prescribed porous-medium application, both schemes resolve
pulse attenuation, immobile storage and delayed release after
injection ceases. Separate refinements of detector traces and mass
outputs support these responses.

Future work should develop sharper a priori spatial estimates and
error bounds for detector traces and captured boundary fluxes.
These estimates could guide adaptive selection of the spatial mesh,
time step and number of exchange compartments. Comparisons at
matched accuracy, using a consistent accounting of computational
cost, are also needed to assess the relative efficiency of the
explicit and implicit schemes.
\section*{Author contributions}
Both authors contributed equally to this work.
\section*{Data availability}
The accompanying package contains the solver code, computed arrays and
full-precision records used in this study. It also includes run
configurations, source hashes and all reported figures.
\section*{Declarations of interest}
None.
\section*{Funding information}
This research did not receive any specific grant from funding agencies
in the public, commercial, or not-for-profit sectors.

\section*{Supplementary material}
The numerical studies can be run with \texttt{run\_paper.py} in fresh mode.
The package documents the separate drivers and their refinement schedules.
Saved mode reproduces figures and tables from computed data.
The combined-bound check evaluates the reported estimates on those trajectories.
\setlength{\bibsep}{2pt}

\end{document}